\input ./style/arxiv-general.cfg
\documentclass[MSNbibl,number,citesort,secthm,seceqn,nameyear,dvips]{arxbj}
\makeatletter
   \@ifpackageloaded{graphicx}{}{\usepackage{graphicx}}
\makeatother
\volume{22}
\issue{3}
\pubyear{2016}
\firstpage{1937}
\lastpage{1961}
\doi{10.3150/15-BEJ715}% Updated by VTEXPTS2LaTeX.exe, 20.05.2015 15:34
\docsubty{FLA}

\makeatletter
\newremark{cond}{RE condition}
\newcommand{\eqref}[1]{(\ref{#1})}
\newremark{rem}{Remark}[section]
\newremark{exa}{Example}
\newtheorem{prop}{Proposition}[section]
\newtheorem{lemma}{Lemma}[section]
\newtheorem{cor}{Corollary}
\makeatother

\begin{document}
\begin{frontmatter}

\title{An analysis of penalized interaction models}
\runtitle{An analysis of penalized interaction models}

\begin{aug}
%%%% inicialai - be tarpu
% Corresponding author: Junlong Zhao - zhaojunlong928@126.com% Updated by VTEXPTS2LaTeX.exe, 21.05.2015 08:13
%by VTEXPTS2LaTeX.exe, 20.05.2015 15:34
\author[A]{\inits{J.}\fnms{Junlong}~\snm{Zhao}\corref{}\thanksref{A}\ead[label=e1]{zhaojunlong928@126.com}}%,
%\author[]{\inits{}\fnms{}~\snm{}\thanksref{}\ead[label=]{}}
 \and
\author[B]{\inits{C.}\fnms{Chenlei}~\snm{Leng}\thanksref{B}\ead[label=e2]{C.Leng@warwick.ac.uk}}
%%\runauthor{} %% auto
%\dedicated{}
\address[A]{School of Mathematics and Systems Science, Beihang
University, LMIB of the Ministry of Education, China. \printead{e1}}
\address[B]{Department of Statistics,
University of Warwick, UK. \printead{e2}}
\end{aug}

% HISTORY:
%
\received{\smonth{4} \syear{2014}}% Updated by VTEXPTS2LaTeX.exe,
%20.05.2015 15:34
%
\revised{\smonth{11} \syear{2014}}% Updated by VTEXPTS2LaTeX.exe,
%20.05.2015 15:34

% ABSTRACT
\begin{abstract}
An important consideration for variable selection in interaction models
is to
design an appropriate penalty that respects hierarchy of the importance
of the
variables. A common theme is to include an interaction term only after the
corresponding main effects are present. In this paper, we study several
recently proposed approaches and present a unified analysis on the convergence
rate for a class of estimators, when the design satisfies the restricted
eigenvalue condition. In particular, we show that with probability
tending to
one, the resulting estimates have a rate of convergence $s\sqrt{\log p_1/n}$
in the $\ell_1$ error, where $p_1$ is the ambient dimension, $s$ is
the true
dimension and $n$ is the sample size. We give a new proof that the restricted
eigenvalue condition holds with high probability, when the variables in the
main effects and the errors follow sub-Gaussian distributions. Under this
setup, the interactions no longer follow Gaussian or sub-Gaussian
distributions even if the main effects follow Gaussian, and thus
existing works are not applicable. This result is of independent interest.
\end{abstract}

% KEYWORDS
% visi is mazosios raides ir pagal abecele
\begin{keyword}
\kwd{convergence rate}
\kwd{hierarchical variable selection}
\kwd{high-dimensionality}
\kwd{interaction models}
\kwd{Lasso}
\kwd{restricted eigenvalue condition}
\end{keyword}
\end{frontmatter}

%s1 #&#
\section{Introduction}\label{sec1}
High-dimensional datasets are predominantly characterized by a large
ambient dimension of the covariates $p$ and a small number of the
observations $n$. An important assumption in analyzing such datasets is
that the true dimension of the relevant covariates $s$ is often smaller
than $n$. Let $(X_i,Y_i), i=1,\ldots,n$ be i.i.d. observations where
$X_i=(X_{i1},\ldots, X_{ip})^T$ is a $p$-dimensional regressor and
$Y_i$ is a scalar response. Denote $\mathbb{X}=(X_1,\ldots,X_n)^T\in
R^{n\times p}, \mathbb{Y}=(Y_1,\ldots, Y_n)^T\in R^{n}$. Consider the usual
linear model
%
%e1.1 #&#
\begin{equation}
\label{LM}
\mathbb{Y}=\mathbb{X}\alpha+\varepsilon,
\end{equation}
where $\alpha\in R^p$ is a sparse vector, $\varepsilon=(\varepsilon
_1,\ldots, \varepsilon_n)^T$ and $\varepsilon_i$ are i.i.d.
errors. There are a number of important approaches for recovering the
sparse vector $\alpha$ from this model, such as
Lasso (Tibshirani \cite{Tibshirani1996}), the SCAD (Fan and Li \cite{FanLi2001}), the Dantzig selector
(Candes and Tao \cite{CandesTao2007}) and many others. The main idea is to make use
of a sparsity-encouraging penalty that gives coefficient estimation and
variable selection at the same time. In particular, the popular Lasso
method aims to estimate $\alpha$ as
\[
\hat{\alpha}=\operatorname{arg}\min_{\alpha} \|\mathbb{Y}-\mathbb {X}
\alpha\| _2^2+\lambda_n\|\alpha
\|_1,
\]
where $\lambda_n$ is a penalty parameter suitably chosen and $\|\cdot
\|_q$ is the $\ell_q$ norm of a vector. A fundamental problem in
analyzing such approaches is to derive sufficient conditions such that
the resulting estimates are consistent with appropriate rate of
convergence.  Zhang and Huang \cite{ZhangHuang2008} analyzed the sparsity and bias of
this formulation. Meinshausen and Yu  \cite{MeinshausenYu2009} showed that with high
probability, all important variables are selected by the Lasso.
In an important work,
Bickel, Ritov and Tsybakov \cite{Bickeletal2009} proposed the restricted eigenvalue (RE)
condition. van~de Geer and B{\"u}hlmann \cite{GeerBuhl2009} considered a weaker compatibility
condition to guarantee the consistency of Lasso in the $\ell_1$ or
$\ell_2$ sense. Subsequently, many authors studied sufficient
conditions under which the restricted eigenvalue condition holds.
Raskutti, Wainwright and Yu
\cite{Raskuttietal2010} proved that the RE condition holds with large
probability for Gaussian designs, and Rudelson and Zhou \cite{Zhou2011} showed that the
RE condition holds for sub-Gaussian distributions with large probability.

In practice, it is often necessary, sometimes even mandatory, to
consider interactions of the main variables and to select the variables
in a hierarchical fashion. For example, in genetic studies, the main
effects of genetic and environmental factors often fail to explain
enough phenotypical variation. Gene-environment interactions are
considered pivotal in this context. In analysis of variance models,
interactions are often analyzed only when the corresponding main
effects are present. In these settings, we consider not only the main
effects of $X_{i1},\ldots, X_{ip}$ but also the effects of their
interactions $X_{ik_1}X_{ik_2}$, for any $k_1, k_2 =1,\ldots, p$ with
$k_1< k_2$. Let $Z_i=(X_i^T, X_i^{*T})^T$ with $X_i=(X_{i1},\ldots,
X_{ip})^T\in R^p$ and $X_i^{*}=(X_{i1}X_{i2},\ldots,X_{ik_1}X_{ik_2},
\ldots,X_{i(p-1)}X_{ip})^T$, $i=1,\ldots,n$ and
$\mathbb{Z}=(Z_1,\ldots
,Z_n)^T\in R^{n\times p_1}$ with $p_1=p(p+1)/2$.
Similar to the model that only contains the main effects, we write an
interaction model as
\[
\mathbb{Y}=\mathbb{Z}\beta+\varepsilon,
\]
where $\beta\in R^{p_1}$ and $\mathbb{Z}\in R^{n\times p_1}$ with columns
consisting of the main effects and interaction terms. In this paper, we
only consider two-way interaction models.
It is unclear if the main results can be extended to higher-order
interaction models.

It is generally desirable to select an interaction effect only after the
corresponding main effects are selected. As Bien, Taylor and Tibshirani \cite{Bienetal2013} and the
references therein argued, an interaction should be allowed into the
model only if the corresponding main effects are present. A version of this
statement can be found in Assumption (A1) in Section~\ref{sec2}.
To achieve this goal, several methods have been developed. Yuan, Joseph and Zou \cite{Yuanetal2009} proposed to use a two-step procedure after obtaining a
first-step estimate that is consistent. The theoretical properties are
studied when $p$ is fixed.
Radchenko and  James \cite{RadchenkoJames2010}
used a nonconvex
penalty for enforcing this strong heredity constraint. A few approaches
formulated as convex optimization were also studied.
Zhao, Rocha and Yu \cite{Zhaoetal2009} proposed to use the CAP approach and
Radchenko and James \cite{RadchenkoJames2010} proposed the VANISH method for hierarchical
selection of variables. Bach \textit{et~al.} \cite{Bachetal2012} suggested a framework
for structured penalty. Bien, Taylor and Tibshirani \cite{Bienetal2013} studied a constrained
optimization model for the same purpose.
These three papers (Zhao, Rocha and Yu \cite{Zhaoetal2009},
Choi, Li and Zhu \cite{Choietal2010},
Bach \textit{et~al.} \cite{Bachetal2012}) consider convex optimization problems of the following type:
\[
\hat{\beta}=\operatorname{arg}\min_{\theta} \|\mathbb{Y}-\mathbb {Z}
\theta\| _2^2+\lambda_n P_e(
\theta),
\]
where $P_e(\theta)$ is a convex penalty usually consisting of norms on
groups of variables Yuan and Lin \cite{YuanLin2006}, Lin and Zhang \cite{LinZhang2006}.
A detailed discussion of this penalty is given in Section~\ref{sec2}.
A desirable property of this one-step formulation is that convex
optimization techniques can be exploited to yield a global optimal solution.
Among these convex formulations, only Radchenko and James \cite{RadchenkoJames2010}
studied theoretical properties of their estimate in terms of model
selection consistency when the design matrix satisfies strong
conditions, similar to the irrepresentable conditions in Zhao and Yu \cite{ZhaoYu2006}.
It is thus of great interest to explore other theoretical aspects of a
selected interaction model.

Although the convex optimization methods are attractive conceptually,
the consistency of the estimates for large $p$ in terms of the
convergence rate is largely unknown, especially under weaker conditions
than those in Radchenko and James \cite{RadchenkoJames2010}. The main goal of this paper
is to
consider this important issue for a class of penalties. Our results
give a general method for proving consistency when interactions are
present. In particular,
two major contributions are made in this paper.

First, assuming sub-Gaussian distributions on the error and the
variables in
the main effects, we establish the consistency of penalized estimators in
interaction models proposed by
Zhao, Rocha and Yu \cite{Zhaoetal2009},
Radchenko and James \cite{RadchenkoJames2010}, Bach \textit{et~al.} \cite{Bachetal2012} and
Bien, Taylor and Tibshirani \cite{Bienetal2013}. In particular, we show that under appropriate
conditions and a modified restricted eigenvalue (RE) condition
Bickel, Ritov and Tsybakov \cite{Bickeletal2009}, these penalized estimates achieve the rate of
convergence in the order $s\sqrt{\log p_1/n}$, which matches the
optimal rate
of the Lasso when no interaction is considered. These results also
motivate us
to develop a unified analysis of a general class of methods for penalized
variable selection in interaction models. The main results are nonasymptotic
and presented in Theorem~\ref{th3}. We note that Negahban \textit{et~al.} \cite{Negahban2012} put
forward a unified theory for high-dimensional $M$-estimation when the penalty
satisfies a notion of decomposability. However, as will become clear
later, the
penalties we investigate in this paper are not necessarily
decomposable, and thus the
theory developed in Negahban \textit{et~al.} \cite{Negahban2012} is not directly applicable to
the interaction
models under consideration.

Second, we generalize the RE condition to the case where both the main
and the interaction terms are presented in covariate $Z_i$.
We remark that the presence of the interactions poses significant
theoretical difficulties. For example, when the main effect covariates
follow a multivariate Gaussian distribution, the joint distribution of
the main effects and the interactions are no longer Gaussian, or even
sub-Gaussian.
Thus, the regularity conditions of some existing results are violated,
for example, in Raskutti, Wainwright and Yu \cite{Raskuttietal2010} and  Rudelson and Zhou \cite{Zhou2011}. A~new
analysis in Theorem~\ref{RENnorm} is provided to address this issue.
This result can be of independent interest.

There is a large literature more recently devoted to selecting
variables in interaction models. One particular class is based on
identifying a superset of the variables that includes the true model as
a subset.  Hao and Zhang \cite{HaoZhang2012a,HaoZhang2012b,HaoZhang2014}
utilized an interesting algorithm to build a hierarchical model using
the forward selection  (Wang  \cite{Wang2009}). Another class is to exploit
algorithms for interaction selection. Shah \cite{Shah2012} exploited a
backtracking iterative algorithm to identify a true model.
Hall and Xue  \cite{HallXue2013} adopted a simple recursive approach to allow
identifications. Bickel, Ritov and Tsybakov  \cite{Bickeletal2010} studied a procedure with a
sequential Lasso fit.

The main content of this paper is arranged as follows. In Section~\ref{sec2}, we
establish the consistency of the estimators with penalty functions
defined in Zhao, Rocha and Yu \cite{Zhaoetal2009}, Radchenko and James \cite{RadchenkoJames2010}, and
Bien, Taylor and Tibshirani \cite{Bienetal2013}. We show that our result is applicable to a
wider class of penalties appropriately defined. In Section~\ref{sec3}, we
generalize the RE condition to the interaction
model setting. We develop sufficient conditions to guarantee the RE
condition to hold. The proofs are delayed to Appendix~\ref{appA} and some
auxiliary results are presented in Appendices~\ref{appB}~and~\ref{appC}.

%s2 #&#
\section{A theory for penalized estimation in interaction
models}\label{sec2}

Recall that $p_1=p(p+1)/2$ is the number of the variables in an
interaction model with $p$ main effect covariates. Consider the model
\[
\mathbb{Y}=\mathbb{Z}\beta+\varepsilon,
\]
where $E(\varepsilon)=0$ and $\mathbb{Z}=(Z_1,\ldots,Z_n)^T$ is an
$n$ by $p_1$ matrix. The columns of $\mathbb{Z}$ consist of the main
effects and the second-order interactions as
$X_{ik_1}X_{ik_2},~k_1<k_2$, and $Z_i$ is independent of $\varepsilon_i$.
Without loss of generality, we assume that $E(X_i)=0$ and
$E(X_{ik_1}X_{ik_2})=0$ for $k_1<k_2$. Otherwise, we can replace $X_i$
and $X_{ik_1}X_{ik_2}$ by $X_i-E(X_i)$ and
$X_{ik_1}X_{ik_2}-E(X_{ik_1}X_{ik_2})$, respectively.
Let $\operatorname{cov}(X_i)=\Sigma_x$ and $\lambda_{\max,x},
\lambda_{\min,x}$
be the corresponding largest and smallest eigenvalues, respectively.
Similarly, we denote $\operatorname{cov}(Z_i)=\Sigma_z$ and define
$\lambda_{\max,z}, \lambda_{\min,z}$.
Denote the true parameter as $\beta=(\beta^{(1)T},\beta^{(2)T})^T$
with $\beta^{(1)}=(\beta_1,\ldots, \beta_p)^T\in R^p$
and $\beta^{(2)}=(\beta_{12}, \ldots, \beta_{(p-1)p})^T\in
R^{p_1-p}$. Denote support sets of $\beta^{(1)}$, $\beta^{(2)}$ and
$\beta$ as $\mathcal{S}^{(1)}$, $\mathcal{S}^{(2)}$ and $S$, respectively.
Then it is clear that $S=\mathcal{S}^{(1)}\cup\mathcal{S}^{(2)}$.
Let $s=|S|$ be the
cardinality of $S$.
To simplify the description, we also denote $\beta=(\beta_1,\ldots
,\beta_p,\ldots,\beta_{p_1})^T$
and use a single index $j$ to denote either a main effect or an interaction.
We focus on the strong hierarchy principle Yuan, Joseph and Zou \cite{Yuanetal2009},
Choi, Li and Zhu \cite{Choietal2010} as stated in the following assumption.
\begin{longlist}[(A1)]
\item[(A1)] If $(i,j)\in\mathcal{S}^{(2)}$ then both $i \in\mathcal
{S}^{(1)}$ and
$j\in\mathcal{S}^{(1)}$.
\end{longlist}
Assumption (A1) means that if an interaction is selected, that is, if
$\beta_{ij}\ne
0$, then the corresponding main effects should also be selected, that
is, $\beta_i\ne0$ and
$\beta_j\ne0$, for any $1\le i,j\le p$. This assumption is very
natural to
achieve hierarchical selection and is a well-accepted practice in
statistics. See, for example, Section~1.2 of Bien, Taylor and Tibshirani \cite{Bienetal2013}.
However, this qualitative constraint also makes variable selection
for interaction models extremely difficult. From a computational point of
view, if one is content with selecting variables via the Lasso penalty as
$P_e(\beta)=\sum_j|\beta_j|$, enforcing Assumption (A1) in this Lasso
formulation would give a combinatorial problem that is challenging to
solve. This is why a convex relaxation of (A1) such as that in
Zhao, Rocha and Yu \cite{Zhaoetal2009} is attractive and tractable. From a theoretical
viewpoint, having a hierarchical structure also complicates the theoretical
analysis. As will be seen later, the convex relaxations of Assumption
(A1) in the literature eventually lead to  {overlapping} group Lasso
penalties, for which the familiar techniques such as those in
Bickel, Ritov and Tsybakov \cite{Bickeletal2009}
and Negahban \textit{et~al.}   \cite{Negahban2012} fail to work.
Finally,
the main idea and analysis in this paper may also be generalized to the
so-called weak hierarchy where $(i,j)\in\mathcal{S}^{(2)}$ leads to
either $i
\in
\mathcal{S}^{(1)}$ or $j\in\mathcal{S}^{(1)}$. Since the majority of
the works in the
literature deal with strong hierarchy, we focus on this case in the
following.

For hierarchical selection of variables that respects the strong
heredity principle, many authors Zhao, Rocha and Yu \cite{Zhaoetal2009}, Radchenko and James \cite{Bienetal2013} considered the
penalized least squares problem
\[
\hat{\beta}=\operatorname{arg}\min_{\theta} \biggl\{\frac
{1}{2n}\|\mathbb{Y}-\mathbb{Z}
\theta\|_2^2+\lambda_n P_e(\theta
) \biggr\},
\]
where $P_e(\theta)$ is a convex penalty.
Let $\theta=(\theta_1,\ldots,\theta_{p},\theta_{12},\ldots,\theta
_{(p-1)p})^T$. The general penalty used by the composite absolute
penalty method (Zhao, Rocha and Yu \cite{Zhaoetal2009}, CAP) entertains
%
%e2.1 #&#
\begin{equation}\label{RJ-pen}
P_e(\theta)=\sum_{j=1}^p \bigl\|(
\theta_j, \theta_{jk \dvt j<k}, \theta _{kj \dvt k < j})
\bigr\|_q + \sum_{j<k} |\theta_{jk}
|,
\end{equation}
where $q>1$. The same penalty is studied in the framework of structured
sparsity by Bach \textit{et~al.} \cite{Bachetal2012}.
A special case is thoroughly investigated
by  Radchenko and James \cite{RadchenkoJames2010} for $q=2$.
Bien, Taylor and Tibshirani \cite{Bienetal2013} proposed a constrained Lasso formulation for
hierarchical selection, which is equivalent to using the penalty
%
%e2.2 #&#
\begin{equation}
\label{Bien}
P_e(\theta)=\sum_{j=1}^p
\max\bigl\{|\theta_j|, \bigl\|(\theta_{jk \dvt j<k},
\theta_{kj \dvt k < j})\bigr\|_1\bigr\}+\sum_{j<k}|
\theta_{jk}|,
\end{equation}
which is similar in nature to but different in form from \eqref{RJ-pen}.

As will be shown later, $P_e(\theta)$ is upper and lower bounded by
$c\|\theta\|_1$ for some constant $c>0$. However, different amounts of
penalty are applied on the terms in (\ref{RJ-pen}) or (\ref{Bien}) to achieve the
effect of hierarchical selection. Thus, the penalties in (\ref{RJ-pen}) and
(\ref{Bien}) are useful for selecting variables hierarchically, while the
Lasso penalty $\|\theta\|_1$ ignores the relative importance of the
main effects and interactions.

We note that Negahban \textit{et~al.} \cite{Negahban2012} developed a unified framework for
$M$-estimator in high-dimensional setting if the regularizer is decomposable
in the sense that $P_e(\theta)=P_e(\theta_A)+P_e(\theta_{A^c})$ for some
index set $A$, where $\theta_A$ and $\theta_{A^c}$ are subvectors of
$\theta$. Decomposable regularizers encompass many useful penalties such
as the $\ell_1$ penalty of Lasso as special cases. As general as the
decomposable assumption is, the penalties in \eqref{RJ-pen} and \eqref{Bien}
are not decomposable expect for some trivial cases. Indeed, it was
pointed out
by Negahban \textit{et~al.} \cite{Negahban2012} that their framework does not apply to the
group Lasso
with overlapping groups, a special case of which corresponds to the
interaction model we study. Therefore, the conclusion of  Negahban \textit{et~al.} \cite{Negahban2012}
is not applicable here and new theory needs to be developed.

To study the consistency of the penalty functions defined in (\ref{RJ-pen}), and (\ref{Bien}), we make the following assumption.
\begin{longlist}[(A2)]
\item[(A2)] There exists fixed constants $a_1,a_2$ such that the
eigenvalues of $\operatorname{cov}(Z_i)$ satisfies $0<a_1<\lambda_{\min
,z}\le\lambda_{\max,z}<a_2<\infty$.
\end{longlist}
Assumption (A2) is made on the population covariance of $\mathbb{Z}$,
and is
generally made in the literature
Wang \cite{Wang2009}, Raskutti, Wainwright and Yu \cite{Raskuttietal2010}.
In particular, (A2) implies that the eigenvalues of $\operatorname{cov}(X_i)$
satisfies $0<\lambda_{\min,x}\le\lambda_{\max,x}<\infty$ and that
the diagonal elements of $\Sigma_z$, denoted by $\{h_j^2, 1\le j\le
p_1\}$ are finite, that is, $\max_j h_j\le h_0<\infty$ for some
$h_0$. From the proof of the theoretical results, Assumption (A2) can
be relaxed to allow $\lambda_{\min,z}$ to tend to zero and $\lambda
_{\max,z}$ to~$\infty$, respectively, at the expense of slower
convergence rates for the estimates. The following Proposition~\ref
{A2fornormal} shows that (A2) holds if $X_i$ is normal.

%prthm.1 #&#
\begin{prop}\label{A2fornormal}
If $X_i$ follows a multivariate normal distribution with $0<\tilde{a}_1<
\lambda_{\min,x}\le\lambda_{\max,x}<\tilde{a}_2<\infty$, for some
fixed $\tilde{a}_1,\tilde{a}_2$, then \textup{(A2)} holds.
\end{prop}

For simplicity, let $v=\hat{\beta}-\beta$ be the difference between
the true $\beta$ vector and the estimate obtained with penalty in
(\ref{Bien}) or (\ref{RJ-pen}). We are interested in the convergence
rate of $v$ when both $p$ and $n$ go to infinity.
Denote $\|a\|_\Sigma=(a^T\Sigma_z a)^{1/2}$ for any $a\in R^{p_1}$.
Let $S^{p_1-1}=\{a \in R^{p_1} \dvt \|a\|_2=1\}$ be the $p_1$-dimensional
sphere. For any vector $u\in R^{p_1}$ and any set $A\subseteq\{
1,2,\ldots,p_1\}$, we denote $u_A \in R^{p_1}$ as the vector with the
$j$th element as $u_j$ if $j\in A$, and 0 otherwise.

We develop the theory under the assumption that $X_i$ and $\varepsilon
_i$, $i=1,\ldots,n$ are i.i.d. observations from sub-Gaussian
distributions. Let
\[
K_x=\|X_i\|_{\psi_2}=\sup_{\|u\|=1}
\bigl\|u^TX_i\bigr\|_{\psi_2}=\sup_{\|u\|
=1}
\sup_{q\ge1}~q^{-1/2}\bigl(\mathrm{E}\bigl|u^TX_i\bigr|^q
\bigr)^{1/q}
\]
be the $\psi_2$ norm of vector $X_i$, and similarly denote
$K_e=\|\varepsilon_i\|_{\psi_2}$ and $K_{xj}=\|X_{ij}\|_{\psi_2},
1\le j\le
p$.  {This definition of a $\psi_2$
norm comes
from Vershynin \cite{Vershynin2011}}. By the basic property of sub-Gaussian
distribution Vershynin \cite{Vershynin2011}, we have $K_{xj}\le K_x$ for any
$1\le j\le p$. A summary of useful results for sub-Gaussian
distributions are listed in the \hyperref[appA]{Appendix}.

For any $0< \delta\le1$, define the event
\[
\mathcal{A}_0=\bigl\{\bigl\|\mathbb{Z}^T\varepsilon/n
\bigr\|_\infty<C_{e,\delta
}\sqrt{\log p_1/n}\bigr\},
\]
where $C_{e,\delta}=[(1+\eta_0)/c]^{1/2}K_eh_0(1+\delta)$ for some
positive constant $c$ and $\eta_0$, is a constant depending only on
$K_x, K_e$ and $\delta$.
Now we show that event $\mathcal{A}_0$ holds with a probability close
to one.

%prthm.2 #&#
\begin{prop}\label{nonnormA0}
Suppose that $X_i$ and $\varepsilon_i$, $i=1,\ldots, n$,
respectively, are i.i.d. sub-Gaussian. When $\log p_1=\mathrm{o}(n^{1/3})$,
for any $\delta\in(0,1)$, as $n>C_1/\delta$ for some constant
$C_1>0$, it holds that $ P(\mathcal{A}_0)\geq q_{1n}\cdot q_{2n}$,
where $q_{1n}$ and $q_{2n}$ defined in \textup{(\ref{q1})} and \textup{(\ref{q2})},
tend to one as $n\rightarrow\infty$.
\end{prop}

The requirement $\log p_1=\mathrm{o}(n^{1/3})$ is due to the fact that,
generally, an interaction term is not sub-Gaussian any more and is
heavy tailed. Therefore, a larger $n$ is required, compared with
existing results for high-dimensional linear models where only the main
effects are considered and the order is typically $n>s\log p$. The
requirement $\log p_1=\mathrm{o}(n^{1/3})$ here may not be optimal.
However, in this paper, we only assume that $X_i$ is a general
sub-Gaussian variable without further assumption on the distribution of
$X_i$. The requirement on $p$ and $n$ can be seen as the price paid on
the rate in exchange for the generality of our result.

We now provide a sufficient condition that guarantees the consistency
of the penalized estimate. Specially, we
consider the restricted eigenvalue (RE) condition introduced by
Bickel, Ritov and Tsybakov \cite{Bickeletal2009}.

\begin{cond*}
Assume
\[
\mathop{\min_{J\subseteq\{1,\ldots,p_1\}}}_{|J|\leq s} \min_{\|\alpha_{J^c}\|_1\leq k_0 \|\alpha_J\|_1}
\frac{\|
\mathbb{Z}\alpha\|_2}{\sqrt{n}\|\alpha_J\|_2}=M(k_0,s) >0.
\]
\end{cond*}

This version of RE condition was introduced by  Bickel, Ritov and Tsybakov \cite{Bickeletal2009}
to prove the consistency of the Lasso and the Dantzig selector in
absence of the interactions.
There are many different conditions of this type in the literature, imposing
different constraints on the design matrix $\mathbb{X}=(X_1,\ldots,X_n)^T$
for establish convergence results.
Noticeably, Candes and Tao \cite{CandesTao2007} introduced the uniform uncertainty
principle condition on sparse eigenvalues.
Zhang and Huang \cite{ZhangHuang2008} introduced the sparse Reisz condition.
The RE condition defined by Bickel, Ritov and Tsybakov \cite{Bickeletal2009} is weaker than the
that of Candes and Tao \cite{CandesTao2007}. It should also be pointed out that
for both the $\ell_1$ and $\ell_2$ losses, strictly weaker conditions
than the RE were given in
van~de Geer \cite{vandeGeer2007} and
Ye and Zhang \cite{YeZhang2010}. With proper parameters,
the conditions in
Ye and Zhang \cite{YeZhang2010}
include that of Bickel, Ritov and Tsybakov \cite{Bickeletal2009} and that of
van~de Geer \cite{vandeGeer2007} as the special case. In addition,
Raskutti, Wainwright and Yu \cite{Raskutti2011}
also defined a similar restrictive eigenvalue condition to obtain the
minimax rate high-dimensional linear regression over $\ell_q$ ball.

Our RE assumption can be seen as an extension of RE condition of
Bickel, Ritov and Tsybakov \cite{Bickeletal2009} to the interaction model. For the main effects model,
Raskutti, Wainwright and Yu \cite{Raskuttietal2010} proved that the RE condition holds under the
Gaussian assumption, and Rudelson and Zhou \cite{Zhou2011} showed that the RE condition
still holds under the sub-Gaussian assumption. Here, we need the RE
condition for variable $\mathbb{Z}$, which contains both the main
effects and
interaction terms and is clearly not sub-Gaussian. We give sufficient
conditions for the RE condition to hold for the general case Section~\ref{sec3}.
Based on the results of Proposition~\ref{nonnormA0} and the RE
condition, we have the following on the consistency of the penalized estimators.

%th2.1 #&#
\begin{thm}\label{th3}
Consider the penalty defined in \textup{(\ref{RJ-pen})} and \textup{(\ref{Bien})} for
hierarchical selection. Suppose that \textup{(A1)}, \textup{(A2)} and the conditions of
Proposition~\ref{nonnormA0} holds and that the RE condition holds with
$k_0=7$. If $\lambda_n \ge C_{e,\delta}\sqrt{\log p_1/n}$, then with
probability $q_{1n}q_{2n}$, which tends to 1 as $n\rightarrow\infty$,
we have
\[
\|v\|_1\leq\lambda_n D(s) \quad\mbox{and}\quad
P_e(v)\le3\lambda_n D(s),
\]
where $D(s)=16s/ M^2(7,s)$.
\end{thm}

Theorem~\ref{th3} shows that the penalty in
Zhao, Rocha and Yu \cite{Zhaoetal2009}, Radchenko and James \cite{RadchenkoJames2010} and
Bien, Taylor and Tibshirani \cite{Bienetal2013} designed for
hierarchical variable selection all lead to consistent estimate in the
$\ell_1$ norm and the norm defined by $P_e(\cdot)$ under suitable
conditions. In particular, the $\ell_1$ norm $v=\hat\beta-\beta$ is
of the order $s\sqrt{\log p_1/n}$, which matches the rate of
convergence for the Lasso when no interaction is considered. Assumption
(A1) is very important for the convergence rate of the estimate;
otherwise, the estimate may get lower convergence order.
To better understand the behavior of penalty functions introduced by
Zhao, Rocha and Yu \cite{Zhaoetal2009},
Radchenko and James  \cite{RadchenkoJames2010} and
Bien, Taylor and Tibshirani \cite{Bienetal2013}, we consider a more general class of penalty functions.
In fact it can be verified that under assumption (A1), the penalty
functions defined in
Zhao, Rocha and Yu \cite{Zhaoetal2009},
Radchenko and James \cite{RadchenkoJames2010} and
Bien, Taylor and Tibshirani  \cite{Bienetal2013} are special cases of
a general class of penalty functions with $L_1=1$ and $L_2=3$
satisfying the following assumption.
\begin{longlist}[(A3)]
\item[(A3)] Suppose that $P_e(0)=0$. For any $\theta_1,\theta_2\in
R^{p_1}$, we have $P_e(\theta_1+\theta_2)\le P_e(\theta
_1)+P_e(\theta_2)$. For any $ \theta\in R^{p_1}$, $ P_e(\theta)\geq
P_e(\theta_S)+L_1\|\theta_{S^c}\|_1$ with $L_1>1/2$ and $P_e(\theta
_S)\leq L_2 \|\theta_S\|_1$.
\end{longlist}
Assumption $P_e(0)=0$ holds for all penalties designed for variable selection
in the literature. By triangular inequality, the requirement
$P_e(\theta_1+\theta_2)\le P_e(\theta_1)+P_e(\theta_2)$ holds for properly
defined norms or seminorms. We can see (A3) as a weakened form of the
decomposable assumption in Negahban \textit{et al.} \cite{Negahban2012}.
We now give a theorem for any penalty function that satisfies (A3).

%pr1.1 #&#
\begin{prop}\label{prop2}
Assume\vspace*{1pt} that \textup{(A2)}, the conditions of Proposition~\ref{nonnormA0} and
the RE condition hold with $k_0= \frac{2 L_2+1}{2L_1-1}$.\vspace*{1pt} Suppose the
penalty $P_e(\theta)$ satisfies \textup{(A3)}. For the estimator $\hat\beta$
associated with penalty\vspace*{1pt} $P_e(\theta)$, we have $\|v\|_1\leq\lambda_n
D(s)$ with probability at least $q_{1n}q_{2n}$, where $D(s)=\frac{
(L_2+L_1)^2s}{(2L_1-1) M^2(k_0,s)}$.
\end{prop}

Following the proof of this proposition, it is not difficult
to show that the rate of convergence is the same when the $\ell_1$
penalty is used without enforcing the hierarchical structure in (A1).
From this perspective,
Assumption (A3) connects interaction selection with the $\ell_1$
penalty and various hierarchical penalty functions $P_e(\theta)$. The
hierarchical penalty functions typically differentiate the main and
interaction terms by penalizing the main effects $\theta_j$ less and
the interactions $\theta_{ij}$ more. At the same time, $P_e(\theta)$
should behave similarly to $\ell_q,q\le1$, that is used for variable
selection. In fact, letting $\theta=\theta_{S^c}$ and $\theta_S=0$
in (A3), we have $P_e(\theta_{S^c})\ge L_1\|\theta_S^c\|_1$. That is,
although the penalty in $P_e(\theta)$ on $\theta_j, j\notin\mathcal{S}
^{(1)}$ are different from that on $\theta_{ij}, (i,j)\notin\mathcal{S}
^{(2)}$, both behave similarly to the $\ell_1$ penalty. On the other
hand, $P_e(\theta_S)\le L_2\|\theta_S\|_1$, indicating that
the penalty in $P_e(\theta)$ on $\theta_S$ is less than or equal to
that in
the $\ell_1$ penalty. Recall that $S=\mathcal{S}^{(1)}\cup\mathcal
{S}^{(2)}$ and
that (A3) does not assume any specific relationship between $\mathcal{S}
^{(1)}$ and $\mathcal{S}^{(2)}$. As $\mathcal{S}^{(1)}$ and $\mathcal
{S}^{(2)}$ satisfy
(A1), the penalty functions defined by
Zhao, Rocha and Yu
\cite{Zhaoetal2009},
Radchenko and James
\cite{RadchenkoJames2010} and
Bien, Taylor and Tibshirani
\cite{Bienetal2013}
are members of the penalty functions satisfying (A3) with $L_1=1$ and
$L_2=3$. Consequently, it is expected
that the convergence rates using these methods are similar to that with the
$\ell_1$ penalty. However, a key advantage of the hierarchical
penalties is
that the chosen model is more likely to have a hierarchical structure,
a parsimony
principle often desirable in practice. We note, however, that
the rate may not be improvable. To see this,
consider a special scenario where the true model is sparse and only contains
$s\ll p$ main effects. Even one is able to fix all the coefficients of the
interactions as zero, the optimal rate of fitting a main effects model
via the
Lasso penalty is $s\sqrt{\log p/n}$ (Bickel, Ritov and Tsybakov \cite{Bickeletal2009}), the same as the
rate $s\sqrt{\log p_1/n}$ obtained in this paper as $p_1 =\mathrm{O}( p^2)$ in an
interaction model. Whether the rate is optimal is a direction for
future research.

Motivated by Proposition~\ref{prop2}, we can define new penalties such
that the rate of convergence result in Theorem~\ref{th3} is achieved.
We now discuss a few concrete examples when the true model obeys the
strong hierarchical property (A1).

\begin{exa}
For $q>1$, define $P_e(\theta)$ as
\[
(p-1)^{-1}\sum_{1\le j<k\le p} \bigl\{\bigl\|(
\theta_j,\theta_{jk}) \bigr\| _q+\bigl\|(
\theta_k,\theta_{jk})\bigr\|_q\bigr\}+\sum
_{1\le j<k\le p}|\theta_{jk}|.
\]
This penalty follows the strong hierarchical principle and satisfies
(A3) with $L_1=1, L_2=1+(p-1)^{-1}\le2$. Thus, Proposition~\ref
{prop2} is applicable.

This penalty can be seen as an alternative to (\ref{RJ-pen}) and (\ref{Bien}) for hierarchical
variable selection, where the factor $p-1$ reflects the fact that every main
effect appears $p-1$ times in the penalty.
In particular, we construct groups
$(\theta_j,\theta_{jk})$ and $(\theta_k,\theta_{jk})$ and impose
$\ell_q, q>1$
penalty to encourage sparsity in these groups. An additional penalty is imposed
on the interaction as in $|\theta_{jk}|$ to encourage greater sparsity of
$\theta_{jk}$.
\end{exa}

\begin{exa}
Due to the connection between variables, for
example, in times series analysis, we want to select variables and
their interactions in a contiguous pattern. The idea of a block convex
structure (Bach \textit{et~al.} \cite{Bachetal2012}) is useful for such a setting. To select
main effects and their interactions in the contiguous pattern of order
$d_0$, we can consider the groups as follows.
Let $T_k=(\theta_{km \dvt m>k},\theta_{mk \dvt m<k})$ and construct groups as
$G_j=(\theta_{j-d_0+1},\ldots, \theta_j, H_j)$ with
$H_j=(T_{j-d_0+1}, \ldots, T_j)$ for $j\ge d_0$.
Then for $q>1$, we define the penalty
\[
P_e(\theta)=\sum_{j\ge d_0}^p \bigl(
\|G_j\|_q+{2^{-1}}\|H_j\|
_1 \bigr).
\]
Clearly, if $d_0=1$, we have the penalty defined in (\ref{RJ-pen}).
For a fixed $d_0$, it is easy to verify that $L_1=1$ and $L_2={3}d_0$.

The consideration of using this penalty comes from time series data
analysis. Let $X=(\mathbf{x}_1,\ldots,\mathbf{x}_p)^T\in R^p$ be the
population
version of $X_i$.
A variable $\mathbf{x}_j$ is considered important variable if both
$\mathbf{x}_j$
and the
variables close to $\mathbf{x}_j$ have large main effects and interaction
effects. Therefore, we construct groups, which involve the main effects and
interaction effects of $d_0$ variables that are contiguous, that is,
$\mathbf{x}_{j-d_0+1},\ldots, \mathbf{x}_j$. Note that $T_k$
consists of the
interaction effect
of $\mathbf{x}_{k}$ with all the other variables. The group $G_j$ is the
collection $d_0$ contiguous
main effect $\theta_{j-d_0+1},\ldots, \theta_j$ and their interactions
$T_{j-d_0+1}, \ldots, T_j$ with all the other variable in the study.
If the
tuning parameter is chosen appropriately, this penalty is more likely
to give
a model with either contiguous sequences of variables with main
effects, or
contiguous sequences of variables with their main effects and
interactions among those nonzero
variables.
\end{exa}

\begin{exa}
Suppose that we want to consider the
interactions of the main effects in a nested manner, that is, we want
to consider the interaction of the $j$th variable $\mathbf{x}_{j}$ and the
$k$-variable $\mathbf{x}_{k}$ with $k< j$ in a contiguous setting. We
can define
$\tilde{G}_j=(\theta_{j}, \tilde{H}_j)$, where $\tilde{H}_j=(\theta
_{kj \dvt k< j})$ for {$1\le j\le p$}.
For $q>1$, we can define
\[
P_e(\theta)={\sum_{j=1}^p} \bigl(\|
\tilde{G}_j\|_q+\|\tilde {H}_j
\|_1 \bigr).
\]
For this penalty satisfying the strong hierarchical principle, it is
easy to verify that (A3) holds with $L_1=1$ and $L_2=2$. Thus, the
conclusion of Proposition~\ref{prop2} and that of Theorem~\ref{th3} hold.

If the tuning parameter is chosen appropriately, this penalty would
give a
model consisting of $\theta_j$ and all (or some) of its interactions
with variables $\mathbf{x}_k$ as
$\theta_{kj}$ for $k<j$. More precisely, this penalty would give
models with main
effects and its interactions with variables that have smaller indices.
This is
interesting when variables are organized along time. When a
variable has no effect, all its interactions with preceding variables have
zero effect but its interactions with later variables may be nonzero.
\end{exa}

%s3 #&#
\section{The RE condition}\label{sec3}

In our setting, variable $Z_i$, containing both the main and
interaction terms of $X_i$, is not sub-Gaussian if $X_i$ follows
multivariate normal. Thus, we need to develop sufficient conditions for
the RE condition to hold. The previous results in the literature
assuming sub-Gaussianity are not applicable here.

To simplify our proof, we first introduce a lemma. Let $X,X^*, Z$ be
the population version of $X_i,X_i^*, Z_i$. We denote $Z-E(Z)=(\tilde
{X}^T,\tilde{Z}^T)^T$, where $\tilde{X}=X-E(X)$ and $\tilde
{Z}=X^{*}-E(X^*)$ and for any $u\in S^{p_1-1}$,
$u=(u^{(1)T},u^{(2)T})^T$ with $u^{(1)}\in R^{p}$ and
$u^{(2)}=(w_{12},\ldots, w_{1p}$, $w_{23},\ldots,  w_{2p}, \ldots,
w_{(p-1)p})^T\in R^{p_1-p}$. Define $w_{ji}=w_{ij}$ for $i<j$ and
$w_{ii}=0$, $i=1,\ldots, p$ and construct the symmetric matrix $W=(w_{ij})$.
It is easy to see that $Z^Tu=X^Tu^{(1)}+2^{-1} X^TWX$.
For any $u\in S^{p_1-1}$, we distinguish three cases: (i) $u^{(1)}=0,
u^{(2)}\ne0$, (ii) $u^{(1)}\ne0, u^{(2)}= 0$ and (iii) $u^{(1)}\ne0,
u^{(2)}\ne0$.
For $u$ in case (iii), we define the correlation coefficient as $\rho
_u=\operatorname{corr}(X^Tu^{(1)}, X^TWX)$, and for that in (i) and (ii),
$\rho_u=0$. Let $X_W^*=[\operatorname{var}(X^TWX)]^{-1/2}X^TWX$ for
case (i) and
(iii), and $X_W^*=0$ for case (ii). We make the following assumptions.
\begin{longlist}[(A5)]
\item[(A4)] We assume $\sup_{u\in S^{p_1-1}} |\rho_u|\le
\bar\rho<1$ for some $\bar\rho$.%, where $\mathcal{E}=\{u: u\in
%S^{p_1-1}\}$.

\item[(A5)] The random vector $X_W^*$ is subexponential with $\psi_1$
norm $K_u$ and there exists $K_0$, such that $\sup_{u\in
S^{p_1-1}} K_u\le K_0<\infty$.
\end{longlist}
It can be verified that (A4) and (A5) hold when $X$ is multivariate
normal. In this case, for any $u\in S^{p_1-1}$, $X^Tu^{(1)}$ is normal and
$X^TWX$ follows a weighted $\chi^2$ distribution; thus the correlation
$\rho_u$ is strictly less than 1. Due to the closeness of the set $S^{p_1-1}$,
we have $\sup_{u\in S^{p_1-1}} |\rho_u|\le\bar{\rho}<1$.
Particularly, as $X\sim N(0,\Sigma_x)$, we have $\rho_u=0$ for any
$u\in S^{p_1-1}$ due to the symmetry.
{In fact, we have shown in the proof of
Proposition~\ref{A2fornormal} that (A4) holds when $X$ is continuous
with $\lambda_{\min,x}>0$.}
Since $X$ is normal, by Magnus and Neudecker \cite{MagnusNeudecker1979}, we have
$\operatorname{var}
(X^TWX)=\operatorname{tr}(W\Sigma_x)^2=\|W\Sigma_x\|_F^2$, where $\|\cdot\|_F$
denotes the Frobenius norm.
Let $W_{\Sigma_x}=\Sigma_x^{1/2}W\Sigma_x^{1/2}$ with associated
eigenvalue decomposition $W_{\Sigma_x}=U \Upsilon U^T$ where $\Upsilon
=\operatorname{diag}(\tau_1,\ldots, \tau_p)$ are the eigenvalues and the
columns of $U$ are eigenvectors. Then $Y=U^T\Sigma_x^{-1/2}X\sim
N(0,I_p)$. Then we have
\[
X_W^*=\|W\Sigma_x\|_F^{-1}
X^TWX=\|W\Sigma_x\|_F^{-1}
Y^T\Upsilon Y =\|W\Sigma_x\|_F^{-1}
\sum_{i=1}^p \tau_i
Y_i^2.
\]
Note that, for $i=1,\ldots,n$, $Y_i^2$'s are independent $\chi^2_1$
that is subexponential with the $\psi_1$ norm $\|Y_i^2\|_{\psi_1}\le
2 \|Y_i\|_{\psi_2}\le2 K_x$, and that $\sum\tau_i^2/\|
W\Sigma_x\|_F^2=1$ for any $u\in S^{p_1-1}$. By
Vershynin \cite{Vershynin2011} (see Proposition~\ref{ApC1} in Appendix~\ref{appC}), for any
$u\in S^{p_1-1}$, $X_W^*$ is subexponential with $K_u=\|Y_i^2\|_{\psi
_1}$. Therefore, it is sufficient to take $K_0=\|Y_i^2\|_{\psi_1}$. We
have the following tail bounds.

%lethm.1 #&#
\begin{lemma}\label{lmRE}
Let $Z_{u,\Sigma}=[Z-E(Z)]^Tu\|u\|_\Sigma^{-1}$. Suppose that
\textup{(A4)}--\textup{(A5)} hold. Then for any $t>0$,
\[
\sup_{u\in S^{p_1-1}} P\bigl(|Z_{u,\Sigma}|>t\bigr)\le C\exp\bigl[- \min
\bigl( t^2/\tilde{K}_0,t/K_0\bigr)\bigr],
\]
where $K_0,\tilde{K}_0$ are constant depend only on $K_x$.
Furthermore, it holds that
\[
\sup_{u\in S^{p_1-1}}
\operatorname{var}\bigl[\bigl(Z^T u\bigr)^{2}\bigr]<\infty.
\]
\end{lemma}

%th3.1 #&#
\begin{thm}\label{RENnorm}
Suppose\vspace*{1pt} that $X$ is sub-Gaussian and the assumptions of
Lemma~\ref{lmRE} hold. For any $0<\epsilon<\lambda_{\min,z}^{1/2}(8\lambda
_{\max,z}^{1/2}+\lambda_{\min,z}^{1/2})^{-1}$, if
\[
\label{nRENnorm}
n>\epsilon^{-1}\max \bigl\{C_1,
\bigl(c_1m_1\tilde{C}_K^{-1}
\log(c_1m_1 p_1)\bigr)^3,
\bigl(4^{-1}\tilde{C}_K^{-1}c_1m_1
\log(4c_1m_1 p_1)\bigr)^3 \bigr
\},
\]
where $m_1=16s(1+k_0)^2$ and $c_1, C_1, \tilde{C}_K$ are positive
constants with $\tilde{C}_K<C_K$,
then the RE condition holds with a high probability that is at least
$\tilde{p}_{0n}$ where $\tilde{p}_{0n}$ defined in \textup{(\ref{tilp0})}
tends to~1, as $n$ goes to infinity.
\end{thm}

When there is no interaction and $X$ is sub-Gaussian, we can show that
the order on the right-hand side of the inequality in Lemma~\ref
{lemm6} can be improved to $\exp(-C_Kn\delta)$. As a result, the
lower bound of the sample size $n$ in Theorem~\ref{RENnorm} can be
improved to a smaller order $\epsilon^{-1}\max\{C_1, c_1m_1\tilde
{C}_K^{-1}\log(c_1m_1 p)$, $ 4^{-1}\tilde{C}_K^{-1}c_1m_1\log
(4c_1m_1 p)\}$. This observation again reflects that selecting
variables in interaction models is more difficult.

%sA #&#
\begin{appendix}

%sA #&#
\section{Proofs of the main results}\label{appA}
We provide the proofs for the theorems. A few additional lemmas that
are needed are given in Appendix~\ref{appB}.

%sA.1 #&#
\subsection{Proof of Proposition \texorpdfstring{\protect\ref{A2fornormal}}{2.1}}
We use the notation defined in Section~\ref{sec3}. Recall that
$Z=(X^T,{X^*}^T)^T$. For any $u\in S^{p_1-1}$,
$u=(u^{(1)T},u^{(2)T})^T$ with $u^{(1)}\in R^{p}$ and
$u^{(2)}=(w_{12},\ldots, w_{1p}$, $w_{23},\ldots,w_{2p}, \ldots,
w_{(p-1)p})^T\in R^{p_1-p}$. Define the symmetric matrix $W=(w_{ij})$
based on $u^{(2)}$ as in Section~\ref{sec3}.
It is easy to see that $Z^Tu=X^Tu^{(1)}+2^{-1} X^TWX$. Then
$\lambda_{\max,z}=\sup_{u\in S^{p_1-1}}\operatorname
{cov}(Z^Tu)$ and
$\lambda_{\min,z}=\inf_{u\in S^{p_1-1}}\operatorname{cov}(Z^Tu)$.

First, we consider $\lambda_{\max,z}$. It holds that
\[
\lambda_{\max,z}\le\sup_{u\in S^{p_1-1}}2 \bigl[\operatorname{var}
\bigl(X^Tu^{(1)}\bigr)+\operatorname{var}
\bigl(2^{-1} X^TWX\bigr)\bigr].
\]
Since $X$ is normal, from Magnus and Neudecker \cite{MagnusNeudecker1979}, we have
$\operatorname{var}
(X^TWX)=\operatorname{tr}(W\Sigma_x)^2=\|W\Sigma_x\|_F^2$, where $\|\cdot\|_F$
denotes the Frobenius norm.
For symmetric positive semidefinite matrices $A$ and $B$, by
Fang,  Loparo and Feng \cite{Fangetal1994}, we have that
$\lambda_{\min,A} \operatorname{tr}(B)\le \operatorname{tr}(AB)\le\lambda_{\max,A} \operatorname{tr}(B)$.
Since $W$ and $\Sigma_{x}$ are symmetric matrices, we have
%
%eA.1 #&#
\begin{equation}
\label{trWSigmax}
\operatorname{tr}\bigl(W^2\bigr)\lambda_{\min,x}^2
\le\operatorname{var}\bigl(X^TWX\bigr)=\operatorname{tr}\bigl(W^2\Sigma
_x^2\bigr)\le \operatorname{tr}\bigl(W^2\bigr)
\lambda_{\max,x}^2.
\end{equation}
Moreover, we have $\operatorname{var}(Z^Tu^{(1)})\le\lambda_{\max
,x}\|u^{(1)}\|_2^2$.
Also, note that $\operatorname{tr}(W^2)+\|u^{(1)}\|_2^2=\|u\|_2^2=1$. Then we have
\begin{eqnarray*}
\lambda_{\max,z} &\le & 2 \sup_{u\in S^{p_1-1}} \bigl[\lambda_{\max
,x}
\bigl\|u^{(1)}\bigr\|_2^2+ 4^{-1}\operatorname{tr}
\bigl(W^2\bigr)\lambda_{\max,x}^2\bigr]\\
&\le & 2\max\bigl
\{ \lambda_{\max,x},\lambda_{\max,x}^2\bigr\}\le2\max\bigl
\{\tilde a_2,\tilde a_2^2\bigr\}.
\end{eqnarray*}

Second,\vspace*{1pt} we consider $\lambda_{\min,z}$. Note that $\lambda_{\min
,x}>\tilde
a_1>0$. Then $X$ is nondegenerate in any direction. For any $u^{(1)}\ne0,
u^{(2)}\ne0$, we have $\operatorname{var}(X^Tu^{(1)})> 0$ and
$\operatorname{var}(X^TWX)> 0$.
Recall the notation $\rho_u=\operatorname{corr}(X^Tu^{(1)}, X^TWX)$ in
Section~\ref{sec3}
which is defined as 0 when $u^{(1)}=0$ or $u^{(2)}=0$.
For any $u\in S^{p_1-1}$ with $u^{(1)}\ne0$ and $u^{(2)}\ne0$, we have
given a simple argument in Section~\ref{sec3} that, when $X$ is normal,
$\sup_{u\in S^{p_1-1}}\rho_u\le\bar\rho<1$, for some $\bar
\rho$. Here,
we show a more general result that if $X$ is continuous with $\lambda
_{\min,x}>0$, this conclusion still holds.

Suppose that on the contrary we have $\sup_{u\in
S^{p_1-1}}\rho_u=1$. Then due to the closeness of $S^{p_1-1}$, we have for
some $u_0=(u_0^{(1)T},u_0^{(2)T})^T$ with $u_0^{(i)}\ne0, i=1,2$, such that
$\rho_{u_0}=1$. Let $W_0=(w_{ij})$ be the matrix constructed from
$u_0^{(2)}$. Then there are some constants $a_0$ and $a_1$ in $R$ such that
$X^Tu_0^{(1)}=a_1(X^TW_0X)+a_0$ holds with probability 1.
Since $u_0^{(1)}\ne0$ and $\lambda_{\min,x}>0$, we have $a_1\ne0$. Since
$W_0$ is a symmetric matrix, its eigenvalue decomposition is denoted as
$W_0=U\Lambda U^T$, where $\lambda=\operatorname{diag}(\lambda_1,\ldots,
\lambda_p)$
is a diagonal matrix of eigenvalues and $U$ consists of the
eigenvectors. Let $u_0^{(1)}=Ub$ for some $b\in R^p$. Then
$V=U^TX:=(\mathbf{v}_1,\ldots,\mathbf{v}_p)^T$ is a nondegenerate
variable in
$R^p$, since $\operatorname{cov}(V)=\operatorname{cov}(X)$ and
$\lambda_{\min,x}>0$.
Therefore, we have with probability\vspace*{-2pt} one,
\[
a_1^{-1}\bigl(V^Tb-a_0
\bigr)=V^T\Lambda V=\sum_{j=1}^p
\lambda_j \mathbf{v}_j^2.
\]
Due to $u^{(2)}\ne0$, we have $W_0\ne0$, and consequently at least one
of $\lambda_i, 1\le i\le p$, is nonzero. Therefore, the right-hand
side is a hyperquadric
and the left-hand side is a hyperplane. Since the set of intersection of
hyperquadric and hyperplane has Lebesgue measure 0 and $V$ is a
continuous and nondegenerate variable in $R^p$, the equation holds with
probability 0. This leads to contradiction.

Furthermore, by the inequality $\sqrt{ab}\le2^{-1}(a+b)$ for any
$a>0,b>0$ and definition of correlation, we have
\begin{eqnarray*}
\operatorname{cov}\bigl(Z^Tu\bigr)&\ge&\operatorname {var}
\bigl(X^Tu^{(1)}\bigr)+\operatorname{var}
\bigl(2^{-1} X^TWX\bigr)-2\bigl|\operatorname{cov}
\bigl(Z^Tu^{(1)}, 2^{-1} X^TWX \bigr)\bigr|
\nonumber
\\[-2pt]
&\ge& (1-\bar\rho)\bigl[\operatorname{var}\bigl(Z^Tu^{(1)}
\bigr)+ \operatorname {var}\bigl(2^{-1} X^TWX\bigr)\bigr]
\nonumber
\\[-2pt]
&\ge& (1-\bar\rho)\bigl[\lambda_{\min,x}\bigl\|u^{(1)}
\bigr\|_2^2+ 4^{-1}\operatorname{tr}\bigl(W^2\bigr)
\lambda_{\min,x}^2\bigr],
\nonumber
%&\ge& (1-\bar\rho)\min(\lambda_{\min,x},\lambda_{\min,x}^2)\nonumber\\
%&\ge&(1-\bar\rho)\min(\tilde a_1,\tilde a_1^2),\nonumber
\end{eqnarray*}
where we use (\ref{trWSigmax}) in the third inequality.
Consequently, it holds that
\[
\lambda_{\min,z}= \inf_{u\in S^{p_1-1}}\operatorname{cov}\bigl(Z^Tu\bigr)
\ge4^{-1}(1-\bar\rho)\min\bigl(\lambda_{\min,x},
\lambda_{\min,x}^2\bigr) \ge4^{-1}(1-\bar\rho)\min
\bigl(\tilde a_1,\tilde a_1^2\bigr),
\]
where we use the fact $\|u^{(1)}\|_2^2+\operatorname{tr}(W^2)=\|u\|_2^2=1$. This
completes the proof.

%sA.2 #&#
\subsection{Proof of Proposition \texorpdfstring{\protect\ref{nonnormA0}}{2.2}}

Recall that $E(Z_{ij})=0$ and $h_j^2=\operatorname{var}(Z_{ij})\le
h^2_0<\infty$,
for $1\le j\le p_1$.
Define $\mathcal{T}= \{\mathbb{Z}\dvt  | (\sqrt{n}\times h_j)^{-1}\|
\mathbb{Z}_j\|
_2-1 |<\delta, \mbox{for all } j=1,\ldots, p_1  \}$. We
now consider the probability of event $\mathcal{T}$.
In fact, let $e_j, j=1,\ldots, p_1$ be the vector with the $j$th
element 1 and 0 elsewhere. Note that for any $1\le j_1<j_2\le p$,
$X_{ij_1}$ and $X_{ij_2}$ are sub-Gaussian with $\psi_2$ norm less
than $K_x$. By Vershynin \cite{Vershynin2011}, $X_{ij_1}X_{ij_2}$ is
subexponential with $\psi_1$ norm no more than $2K_x^2$ (see also
Appendix~\ref{appC}).
By the results of Lemma~\ref{lemm6}, as $n>C_1/\delta$, we have
%
%eA.2 #&#
\begin{eqnarray}
P(\mathcal{T})&=&1-P\bigl(\mathcal{T}^c\bigr)
\nonumber
\\
&\geq& 1-\sum_{j=1}^{p_1} P\bigl(\bigl|
n^{-1}h_j^{-2}\bigl\|\mathbb{Z}e_j
\bigr\|^2_2 -1\bigr|\geq\delta\bigr)
\nonumber
\\[-8pt]
\label{q1}
\\[-8pt]
\nonumber
&\geq& 1- 2p_1\exp \bigl(-C_K (n\delta)^{1/3}
\bigr)
\\
&=&1- 2 \exp \bigl[-C_K(n\delta)^{1/3}+\log
p_1 \bigr]:=q_{1n},
\nonumber
\end{eqnarray}
where $C_K$ depends only on $K_x$ and $h_0$. %
It is clear that $q_{1n}\rightarrow1$, when $n\rightarrow\infty$ due
to the assumption $\log p_1=\mathrm{o}(n^{1/3})$.
On the other hand, letting $V_j=\mathbb{Z}_j^T\varepsilon$,
$j=1,\ldots
,p_1$, we have
\[
P\bigl(\bigl\|\mathbb{Z}^T\varepsilon/n\bigr\|_\infty>t|\mathcal{T}
\bigr)=P\Bigl(\max_j |V_j|>nt|\mathcal{T}\Bigr)
\leq p_1 \max_j P\bigl(|V_j|>nt|
\mathcal{T}\bigr).
\]
By the independence of $Z_i$ and $\varepsilon_i$, we have that
$\mathbb{Z}
_j^T\varepsilon|\mathbb{Z}_j$ is also sub-Gaussian with $\psi_2$ norm
$\tilde{K}_{j,e}=C\|\mathbb{Z}_j\| K_e$ for some constant $C>0$.
Then we have
\[
\max_j P\bigl(|V_j|>nt|\mathbb{Z}_j\bigr)
\leq\max_j e\exp \biggl(-\frac{c
n^2t^2}{K_e^2\|\mathbb{Z}_j\|_2^2} \biggr),
\]
where $c>0$ is constant.
Moreover, we have
\[
\max_j P\bigl(|V_j|>nt|\mathcal{T}\bigr)\leq e \exp
\biggl(-\frac{c
t^2n}{K_e^2(1+\delta)^2h_0^2}\biggr).
\]
Consequently, we have
\[
P\bigl(\bigl\|\mathbb{Z}^T\varepsilon/n\bigr\|_\infty>t|\mathcal{T}
\bigr)\leq p_1e \exp \biggl(-\frac{c t^2n}{K_e^2(1+\delta)^2h_0^2}\biggr).
\]
Taking $t=C_{e,\delta}\sqrt{\log p_1/n}$, with $C_{e,\delta
}=[(1+\eta_0)/c]^{1/2}K_eh_0(1+\delta)$ for any $\eta_0>0$, we have
%
%eA.3 #&#
\begin{equation}
\label{q2}
P\bigl(\bigl\|\mathbb{Z}^T\varepsilon/n
\bigr\|_\infty<C_{e,\delta}\sqrt{\log p_1/n}|\mathcal{T}
\bigr)\geq1-p_1^{-\eta_0}:=q_{2n},
\end{equation}
where $q_{2n}\rightarrow1$ with $n\rightarrow\infty$. Therefore,\vspace*{-2pt} it
holds that
\[
P\bigl(\bigl\|\mathbb{Z}^T\varepsilon/n\bigr\|_\infty\leq
C_{e,\delta}\sqrt{\log p_1/n} \bigr)>q_{1n}q_{2n}.
\]
The proof is completed.

%sA.3 #&#
\subsection{Proof of Theorem \texorpdfstring{\protect\ref{th3}}{2.1}}
\begin{longlist}[\textit{Step} 1.]
\item[\textit{Step} 1.] For the penalty function defined in (\ref{RJ-pen})
and (\ref{Bien}), we show that under (A1), assumption (A3) holds with
$L_1=1$ and $L_2=3$. We take the penalty in (\ref{RJ-pen}) as example
for illustration. The proof for penalty in (\ref{Bien}) can be checked
similarly.

For any $j \notin\mathcal{S}^{(1)}$, by Assumption (A1), pair
$(j,k)\notin
\mathcal{S}^{(2)}$ for all $ j\ne k$.\vspace*{-2pt}
Consequently,
%
%eA.4 #&#
\begin{equation}
\label{eqpr1}
\biggl[|\theta_j|^q+\sum_{k \dvt j<k}| \theta_{jk}|^q+\sum
_{k \dvt k<j}| \theta_{kj} |^q
\biggr]^{1/q}+\sum
_{k \dvt j<
k}|\theta_{jk}|\ge|
\theta_j|+\sum
_{ k \dvt j<k}|\theta_{jk}|.
\end{equation}
And for any\vspace*{-4pt} $j\in\mathcal{S}^{(1)}$,
%
%eA.5 #&#
\begin{eqnarray}
 && \biggl[|\theta_j|^q+\sum_{k \dvt j<k}| \theta_{jk}|^q+\sum
_{k \dvt k<j}| \theta_{kj} |^q
\biggr]^{1/q}+\sum_{k\dvt j<k}|
\theta_{ij}|\nonumber\\[-2pt]
&&\label{eqpr2}\quad\ge \biggl[|\theta_j|^q+ \mathop{
\sum_{k \dvt (j,k)\in\mathcal{S}^{(2)}}}_{j<k} |\theta_{jk}|^q+
\mathop{\sum_{k \dvt (j,k)\in\mathcal{S}^{(2)}}}_{k<j} |
\theta_{kj}|^q \biggr]^{1/q}
\\[-2pt]
\nonumber
&&\qquad{}+\mathop{\sum_{k \dvt (j,k)\in\mathcal{S}^{(2)}}}_{j<k} |\theta
_{jk}|+\mathop{\sum_{k \dvt (j,k)\notin\mathcal{S}^{(2)}}}_{j<k}
|\theta_{jk}|.
\end{eqnarray}
It is easy to see that sum of the left-hand side in (\ref{eqpr1}) and
(\ref{eqpr2}) equals $P_e(\theta)$ and that of the right-hand sides
equals $P_e(\theta_{S})+\|\theta_{S^c}\|_1$. Therefore, $P_e(\theta
)\ge P_e(\theta_S)+\|\theta_{S^c}\|_1$, that is $L_1=1$. Noting the
fact that for $q>1$, $\|\alpha\|_q\le\|\alpha\|_1$ for any vector
$\alpha$, we can see
\begin{eqnarray*}
P_e(\theta)&=&\sum_j \biggl[|
\theta_j|^q+\sum
_{j<k}|
\theta_{jk}|^q+\sum
_{k<j}|
\theta_{kj} |^q \biggr]^{1/q}+\sum_{j<k}|\theta_{jk}|
\nonumber
\\[-2pt]
&\le& \sum_j \biggl[|\theta_j|+\sum_{j<k}| \theta _{jk}|+\sum_{k<j}|
\theta_{kj} | \biggr]+\sum_{
j<k}|
\theta_{jk}|
\nonumber
\\[-2pt]
&=& \sum_{j}|\theta_j|+3\sum_{ j<k}|
\theta _{jk}|
\nonumber
\le3\|\theta\|_1.
\end{eqnarray*}
Therefore, $L_2=3$.

\item[\textit{Step} 2.] Note that the penalty function defined in (\ref
{RJ-pen}) and (\ref{Bien}) are special cases of the function class
satisfying (A3). Consequently, under (A1), (A2) and the RE condition
with $k_0=(2L_2+1)/(2L_1-1)=7$, by Proposition~\ref{prop2}, we have $\|
v\|_1\leq\lambda_n D(s)$ with probability $p_nq_n$, tending to 1 as
$n\rightarrow\infty$. In addition, due to $P_e(\theta)\le3\|\theta
\|_1$ in Step 1, we have $P_e(v)\leq3\lambda_n D(s)$.
This completes the proof.
\end{longlist}

%sA.4 #&#
\subsection{Proof of Proposition \texorpdfstring{\protect\ref{prop2}}{2.3}}
\begin{longlist}[\textit{Step} 1.]
\item[\textit{Step} 1.] By the definition of $\hat{\beta}$ and convexity of $P_e(\theta)$,
we have
\[
\frac{1}{2n} \|\mathbb{Y}-\mathbb{Z}\hat{\beta}\|_2^2+
\lambda_n P_e(\hat{\beta }) \leq\frac{1}{2n} \|
\mathbb{Y}-\mathbb{Z}\beta\|_2^2+\lambda_n
P_e(\beta).
\]
Conditioning on the set $\mathcal{A}_0$, taking $\lambda_n>
2C_{e,\delta
}\sqrt{\log p_1/n}$, and noting that $\hat{\beta}=\hat{\beta
}_S+v_{S^c} $, we have
\begin{eqnarray*}
\frac{1}{n} \|\mathbb{Z}v\|_2^2&\leq& 2
\lambda_n P_e(\beta) -2\lambda_n
P_e(\hat{\beta})+2\frac{v^T\mathbb{Z}^T\varepsilon}{n}
\\
&\leq& 2\lambda_n P_e(\beta) -2\lambda_n
P_e(\hat{\beta})+2\biggl\| \frac{\mathbb{Z}^T\varepsilon}{n}\biggr\|_\infty\|v
\|_1
\\
&\leq& 2\lambda_n \bigl[P_e(\beta) -P_e(
\hat{\beta})\bigr]+ \lambda_n \|v\|_1.
\end{eqnarray*}
And consequently,
%
%eA.6 #&#
\begin{eqnarray}
\frac{1}{n} \|\mathbb{Z}v\|_2^2 &
\leq & 2\lambda_n \bigl[P_e(\beta )-P_e(\hat{
\beta}_S)-L_1\|v_{S^c}\|_1 \bigr]
+\lambda_n \|v\|_1
\nonumber
\\
&\leq& 2\lambda_n \bigl[P_e(-v_S)-L_1
\|v_{S^c}\|_1 \bigr] +\lambda_n \|v
\|_1
\nonumber
\\[-8pt]
\label{Meq1}
\\[-8pt]
\nonumber
&\leq& 2\lambda_n [L_2\|v_S
\|_1-L_1\|v_{S^c}\|_1 ] +\lambda
_n \bigl(\|v_S\|_1+\|v_{S^c}
\|_1\bigr)
\\
\nonumber
& =& \lambda_n \bigl[(2 L_2+1 ) \|v_S
\|_1- (2L_1-1 ) \|v_{S^c}\|_1
\bigr].
\end{eqnarray}
Thus, we have
\[
\|v_{S^c}\|_1\leq\frac{2 L_2+1}{2L_1-1} \|v_S
\|_1:=k_0 \|v_S\|_1.
\]
Recall that the RE condition holds with $M(k_0,s)>0$. Conditioning on
$\mathcal{A}_0$ and this condition, by~(\ref{Meq1}), we have
\begin{eqnarray*}
\|\mathbb{Z}v\|_2^2/n +\lambda_n(2
L_1-1 ) \|v\|_1 &\leq& 2\lambda_n
(L_2+L_1)\|v_S\|_1
\nonumber
\\
&\leq& 2\lambda_n (L_2+L_1)\sqrt{s}
\|v_S\|_2
\nonumber
\\
&\leq& 2\lambda_n (L_2+L_1)\sqrt{s}
\frac{\|\mathbb{Z}v\|_2/\sqrt
{n}}{M(k_0,s)}
\nonumber
\\
&\leq& \frac{\lambda_n^2 (L_2+L_1)^2s}{M^2(k_0,s)} + \biggl( \frac{\|
\mathbb{Z}v\|_2}{\sqrt{n}} \biggr)^2.
\end{eqnarray*}
Therefore, we have
\[
\lambda_n(2 L_1-1 ) \|v\|_1 \leq
\frac{\lambda_n^2 (L_2+L_1)^2s}{M^2(k_0,s)},
\]
and
\[
\|v\|_1 \leq\lambda_n \frac{ (L_2+L_1)^2s}{(2L_1-1)
M^2(k_0,s)}:=
\lambda_n D(s).
\]
\end{longlist}

%sA.5 #&#
\subsection{Proof of Lemma \texorpdfstring{\protect\ref{lmRE}}{3.1}}

The result $\sup_{u\in S^{p_1-1}} \operatorname
{var}[(Z_i^Tu)^{2}]<\infty
$ follows easily from the inequality in Lemma~\ref{lmRE} and
assumption $\lambda_{\max,z}<\infty$ in (A2). We only need to prove
the main inequality in the lemma.
\begin{longlist}[Step 1.]
\item[\textit{Step} 1.]   We have $[Z-E(Z)]^Tu=\tilde{X}^{T}u^{(1)}+\tilde
{Z}^{T}u^{(2)}$. Then by the definition of $W$, it is easy to see
$\tilde{Z}^{T}u^{(2)}=2^{-1} [X^TWX-E(X^TWX)]$.
Then
\begin{eqnarray}
D_u&:=&P\bigl(|Z_{u,\Sigma}|>t\bigr)
\nonumber
\\
&=&P\bigl(\|u\|^{-1}_\Sigma\bigl|\tilde{X}^Tu^{(1)}+
\tilde {Z}^Tu^{(2)}\bigr|>t\bigr)
\nonumber
\\
&\le& P\bigl(\bigl|\tilde{X}^Tu^{(1)}\bigr| > t\cdot2^{-1}
\|u\|_\Sigma\bigr)+P\bigl( \bigl|\tilde {Z}^Tu^{(2)}\bigr|> t
\cdot2^{-1}\|u\|_\Sigma\bigr)
\nonumber
\\
&:=&D_{1,u}+D_{2,u}.
\nonumber
\end{eqnarray}
If $u^{(1)}=0$, it is clear that $D_{1,u}=0$. If $u^{(1)}\ne0$, due to
$\|u^{(1)}\|_2\le\|u\|_2$, it follows that
\begin{eqnarray}
D_{1,u}&=& P\bigl(\bigl|\tilde{X}^Tu^{(1)}\bigr|
\bigl\|u^{(1)}\bigr\|_2^{-1} >t\cdot2^{-1}\|u\|
_\Sigma\bigl\|u^{(1)}\bigr\|_2^{-1}\bigr)
\nonumber
\\
&\le& P\bigl(\bigl|\tilde{X}^Tu^{(1)}\bigr|\bigl\|u^{(1)}
\bigr\|_2^{-1} >t\cdot2^{-1}\|u\| _\Sigma\|u
\|_2^{-1}\bigr)
\nonumber
\\
&\le& P\bigl(\bigl|\tilde{X}^Tu^{(1)}\bigr|\bigl\|u^{(1)}
\bigr\|_2^{-1} > t\cdot2^{-1} \lambda^{1/2}_{\min,z}
\bigr).
\nonumber
\end{eqnarray}
Note that $\tilde{X}$ is sub-Gaussian with $\psi_2$ norm $K_x<\infty
$. Then for any $u^{(1)}\ne0$, $|X^Tu^{(1)}|\|u^{(1)}\|_2^{-1}$ is
sub-Gaussian with $\psi_2$ norm less than $K_x$. Therefore, since
$\lambda_{\min,z}$ is finite, we have
\[
\sup_{u\in S^{p_1-1}} D_{1,u}\le\exp\bigl(- t^2/
\tilde{K}_x^2 \bigr),
\]
where $\tilde{K}_x$ is a constant depending only on $K_x$.

\item[\textit{Step} 2.]  We show the exponential rate of $\sup_{u
\in S^{p_1-1}} D_{2,u}$ and the final conclusion.
\begin{longlist}[\textit{Step} 2.1.]
\item[\textit{Step} 2.1.]  We first show that $\sup_{u\in S^{p_1-1}}
D_W\le L_0<\infty$, where $D_W=[\operatorname{var}(X^TWX)]^{1/2}$. If
$u^{(2)}=0$,
$\operatorname{var}(X^TWX)=0$. We only need to consider $u$ satisfying\vspace*{1pt}
$u^{(2)}\ne0$.
Recall that $\|u\|_\Sigma=(u^T\Sigma_z u)^{1/2}$. It follows that
$ \|u\|_\Sigma^2=u^T\Sigma_z u=\operatorname{var}(X^T
u^{(1)}+2^{-1}X^TWX)\le
\lambda_{\max,z}<\infty$.
That is,
\[
\lambda_{\max,z}\ge\operatorname{var}\bigl(X^T
u^{(1)}+2^{-1}X^TWX\bigr)=A_u^2+B_u^2+2
\rho _u A_uB_u,
\]
where $A_u^2= \operatorname{var}(X^T u^{(1)})$ and
$B_u^2=\operatorname{var}(2^{-1}X^TWX)$. It is
sufficient to show that $\sup_{u\in S^{p_1-1}} B_u$ is bounded.
Suppose that, on the contrary, there exists series $\{u_n, n\ge1\}$,
such that $B_{u_n}\rightarrow\infty$; for simplicity we denote
$B_{u_n}$ as $B_n$, similarly we define $A_n,\rho_n$.
It is easy to see
%
%eA.7 #&#
\begin{equation}
\label{eq4} \lambda_{\max,z}\ge A_n^2+B_n^2+2
\rho_n A_nB_n\ge2(1+\rho_{n})
A_nB_n\ge2(1-\bar\rho)A_nB_n.
\end{equation}
Consequently, $A_n\rightarrow0$ and $\sup_n|\rho_nA_nB_n|\le
\lambda_{\max,z}\bar{\rho}/[2(1-\bar\rho)]$ due to (A4). Taking
$n\rightarrow\infty$, we have $A_n^2+B_n^2+2\rho_n A_nB_n\rightarrow
\infty$. This contradicts with the first inequality in (\ref{eq4}).
Consequently, $\sup_{u\in S^{p_1-1}} B_u$ is bounded.

\item[\textit{Step} 2.2.] First, for any $u\in S^{p_1-1}$, we have $\|u\|
_\Sigma\ge\lambda^{1/2}_{\min,z}>0$.
Noting that $\tilde{Z}^Tu^{(2)}=2^{-1} (X^TWX-E(X^TWX))$, we have
%
%eA.8 #&#
\begin{eqnarray}
D_{2,u}&=& P\bigl(D_W^{-1}\bigl|
\tilde{Z}^Tu^{(2)}\bigr| >t\cdot2^{-1}\|u
\|_\Sigma D_W^{-1}\bigr)
\nonumber
\\
&=&P\bigl(\bigl|X_W^*-E\bigl(X_W^*\bigr)\bigr|>t \|u
\|_\Sigma D_W^{-1}\bigr)
\\
&\le&P\bigl(\bigl|X_W^*-E\bigl(X_W^*\bigr)\bigr|>t\cdot
\lambda^{1/2}_{\min,z} L_0^{-1}\bigr).\nonumber
\end{eqnarray}
Since $X_W^*$ is subexponential with $K_u\le K_0$, by
Vershynin \cite{Vershynin2011}, we have
\[
\sup_{u\in S^{p_1-1}}D_{2,u}\le2 \exp \bigl[-c\min
\bigl(t^2/K_0^2, t/K_0\bigr)
\bigr],
\]
where $c$ is an absolute constant depending on $L_0$ and $\lambda
_{\min,z}$.
We now prove the final conclusion. From Steps 1 and 2, we easily have
\begin{eqnarray}
\sup_{u\in S^{p_1-1}}D_u &\le& \exp\bigl(-t^2/
\tilde{K}_x^2 \bigr)+2 \exp \bigl[-c\min
\bigl(t^2/K_0^2, t/K_0\bigr)
\bigr]
\nonumber
\\
&\le& C\exp\bigl[- \min\bigl( t^2/\tilde{K}_0,t/K_0
\bigr)\bigr],
\nonumber
\end{eqnarray}
where $\tilde{K}_0=\max(\tilde{K}_x,K_0)$.
\end{longlist}
\end{longlist}

%sA.6 #&#
\subsection{Proof of Theorem \texorpdfstring{\protect\ref{RENnorm}}{3.1}}

\begin{longlist}[\textit{Step}~1.]
\item[\textit{Step} 1.]  We first prove that the RE condition holds under the
event $\mathcal{A}_1\cap\mathcal{A}_2$, where $\mathcal{A}_1$ and
$\mathcal{A}_2$ are defined below.
Define the event
\[
\mathbb{B}_1=\bigl\{x \dvt  x\in R^{p_1}, \exists T,
\mbox{ with } |T|\le s, \mbox{ such that } \|x_{T^c}\|_1
\leq k_0 \|x_{T}\|_1 \bigr\}\cap
S^{p_1-1},
\]
where $p_1=p(p+1)/2$ and $S^{p_1-1}$ is the $p_1$-dimensional unit sphere.

Due to the fact that
$\frac{\|\mathbb{Z}x\|_2}{\sqrt{n}\|x_S\|_2} \geq\frac{\|\mathbb
{Z}x\|_2}{\sqrt
{n}\|x\|_2} $, for any $S\subseteq\{1,\ldots, p\}$,
we only need to show that
\[
\min_{|J|\leq s} \min_{\|x_{J^c}\|_1\leq k_0 \|x_J\|
_1} \frac{\|\mathbb{Z}x\|_2}{\sqrt{n}\|x\|_2} >0.
\]
This is equivalent to show
\[
\min_{x\in\mathbb{B}_1} \frac{\|\mathbb{Z}x\|_2}{\sqrt{n}}>0.
\]
Based on the results of Lemma~\ref{lemm5} and $\Lambda_1, \Lambda_2$
defined there, we have $\mathbb{B}_1\subseteq4\operatorname{conv}\Lambda
_1$, where $4 \operatorname{conv}\Lambda_1=\{4x \dvt  x\in\operatorname{conv}\Lambda
_1\}$ and $\operatorname{conv}\Lambda_1$ for the finite set $\Lambda_1$
denotes the convex hull of set $\Lambda_1$.
Consider the events
\[
\mathcal{A}_1= \biggl\{ \biggl| \frac{\|\mathbb{Z}x_0\|_2}{\sqrt{n}\|
x_0\|_\Sigma} -1 \biggr|\leq\epsilon
\mbox{ for all } x_0\in\Lambda_1 \biggr\},
\]
and
\[
\mathcal{A}_2= \biggl\{ \biggl| \frac{\|\mathbb{Z}x\|_2}{\sqrt{n}\|x\|
_\Sigma} -1 \biggr|\leq\epsilon\mbox{ for all } x\in\Lambda_2 \biggr\},
\]
where $\|x_0\|_\Sigma=(x_0^T\Sigma x_0)^{1/2}$ and $\Lambda
_1\subseteq S^{p_1-1}, \Lambda_2\subseteq B_2^{p_1}$ are defined in
Lemma~\ref{lemm5}.
Conditioning on $\mathcal{A}_2$, we have
\[
\frac{\|\mathbb{Z}x\|_2}{\sqrt{n}\|x\|_2} \leq\frac{\|x\|_\Sigma
}{\|x\|
_2}(1+\epsilon) \leq\lambda^{1/2}_{\max,z}(1+
\epsilon) \qquad\mbox{for any } x\in\Lambda_2.
\]
Consequently,
%
%eA.9 #&#
\begin{equation}
\label{eq0} \|\mathbb{Z}x\|_2/\sqrt{n} \leq\lambda^{1/2}_{\max,z}(1+
\epsilon) \|x\| _2 \qquad\mbox{for any }x\in
\Lambda_2.
\end{equation}
By the results of Lemma~\ref{lemm5}, for any $x\in\mathbb{B}_1$,
there exists $x_0\in\Lambda_1$ such that $\|x-x_0\|_2<\epsilon$.
Therefore, we have
%
%eA.10 #&#
\begin{equation}
\label{eq1}
\hspace*{-4pt}\|\mathbb{Z}x_0\|_2/\sqrt{n}-\bigl\|
\mathbb{Z}(x-x_0)\bigr\|_2/\sqrt{n} \leq \|\mathbb{Z}x\|
_2/\sqrt{n} \leq\|\mathbb{Z}x_0\|_2/
\sqrt{n}+\bigl\|\mathbb{Z}(x-x_0)\bigr\| _2/\sqrt{n}.
\end{equation}
Note that $x-x_0\in(\mathbb{B}_1-\mathbb{B}_1)\cap\epsilon
B_2^{p_1}$ and recall from Lemma~\ref{lemm5}
that $ (\mathbb{B}_1-\mathbb{B}_1)\cap\epsilon B_2^{p_1}
\subseteq4\operatorname{conv}(\Lambda_2)$.
Combining with the definition of $\Lambda_2$ and (\ref{eq0}), we have
%
%eA.11 #&#
\begin{equation}
\label{eq2}
\bigl\| \mathbb{Z}(x-x_0)\bigr\|_2/\sqrt{n} \leq4
\sup_{z\in\Lambda
_2} \| \mathbb{Z}z\|_2/\sqrt{n} \leq4
\lambda^{1/2}_{\max,z}(1+\epsilon) \|z\| _2 \leq4
\lambda^{1/2}_{\max,z}(1+\epsilon) \epsilon.
\end{equation}
In addition, conditioning on $\mathcal{A}_1$, we have
%
%eA.12 #&#
\begin{equation}
\label{eq3}
\|x_0\|_\Sigma(1-\epsilon) <\|
\mathbb{Z}x_0\|_2/\sqrt{n} \leq\| x_0\|
_\Sigma(1+\epsilon).
\end{equation}
There, conditioning on $\mathcal{A}_1\cap\mathcal{A}_2$ and
combining (\ref
{eq1})--(\ref{eq3}), as $0<\epsilon<\lambda_{\min,z}^{1/2}(8\lambda
_{\max,z}^{1/2}+\lambda_{\min,z}^{1/2})^{-1}$,
we have, for any $x\in\mathbb{B}_1$,
\[
\|\mathbb{Z}x\|_2/\sqrt{n}\geq\|x_0
\|_\Sigma(1-\epsilon)- 4\lambda ^{1/2}_{\max,z}(1+
\epsilon)\epsilon\ge \lambda^{1/2}_{\min,z}(1-\epsilon)- 4
\lambda^{1/2}_{\max
,z}(1+\epsilon) \epsilon> 0.
\]
That is, with probability $\tilde{p}_{0n}:=P(\mathcal{A}_1\cap
\mathcal{A}_2)$, the
RE condition holds with $M(k_0, s)>0$.

\item[\textit{Step} 2.]  Finally, we show the bound of $\tilde{p}_{0n}$.
In fact, for any $x_0\in S^{p_1-1}$, it holds that
\[
\biggl| \frac{\|\mathbb{Z}x_0\|_2}{\sqrt{n}\|x_0\|_\Sigma}-1 \biggr|\leq \biggl| \frac{\|\mathbb{Z}x_0\|^2_2}{n\|x_0\|^2_\Sigma}-1 \biggr|.
\]
In addition, it is easy to see that
%
%eA.13 #&#
\begin{equation}
\label{A1A2Nnorm}
P(\mathcal{A}_1\cap\mathcal{A}_2)\ge1-
\sum_{x_0\in\Lambda
_1} P \biggl( \biggl| \frac{\|\mathbb{Z}x_0\|^2_2}{n\|x_0\|^2_\Sigma}-1 \biggr|>\epsilon
\biggr)-\sum_{x\in\Lambda_2} P \biggl( \biggl| \frac{\|
\mathbb{Z}x\|
^2_2}{n\|x\|^2_\Sigma}-1 \biggr|\geq
\epsilon \biggr).
\end{equation}
Recall the assumption $E(Z)=0$ and $\Lambda_1\subset S^{p_1-1}$. For
any $x_0\in\Lambda_1$, by Lemma~\ref{lmRE} and consequently by
Lemma~\ref{lemm6}, for any $n>C_1\epsilon^{-1}$ for some constant
$C_1>0$ defined in Lemma~\ref{lemm6},
we have
\[
P \biggl( \biggl| \frac{\|\mathbb{Z}x_0\|^2_2}{n\|x_0\|^2_\Sigma}-1 \biggr|>\epsilon \biggr) \le C_0 \exp
\bigl(-C_K(n\epsilon)^{1/3}\bigr).
\]
Furthermore, recall that $\Lambda_2\subset B_2^{p_1}$. For any $x\in
\Lambda_2$, by letting $\tilde{x}=x\|x\|_2^{-1}\in S^{p_1-1}$, we
note that
\[
\frac{\|\mathbb{Z}x\|^2_2}{n\|x\|^2_\Sigma}=\frac{\|\mathbb{Z}x\|
x\|^{-1}_2\|
^2_2}{n\|x\|x\|_2^{-1}\|^2_\Sigma}=\frac{\|\mathbb{Z}\tilde{x}\|
^2_2}{n\|
\tilde{x}\|^2_\Sigma}.
\]
By the same argument, for any $x\in\Lambda_2$, as $n>C_1\epsilon
^{-1}$ where constant $C_1$ is defined in Lemma~\ref{lemm6}, we have
%
%eA.14 #&#
\begin{equation}
\label{a4}
P \biggl( \biggl| \frac{\|\mathbb{Z}x\|^2_2}{n\|x\|^2_\Sigma}-1 \biggr|>\epsilon \biggr) \le
C_0 \exp\bigl(-C_K(n\epsilon)^{1/3}\bigr).
\end{equation}
Combining (\ref{A1A2Nnorm})--(\ref{a4}), (\ref{cover1}) and (\ref
{cover2}), we have
\begin{eqnarray*}
P(\mathcal{A}_1\cap
\mathcal{A}_2)&\geq & 1-\bigl(|\Lambda_1|+|\Lambda
_2|\bigr)C_0 \exp \bigl(-C_K(n
\epsilon)^{1/3}\bigr)
\\
&\geq & 1-2C_0 \max\bigl\{|\Lambda_1|, |\Lambda_2|
\bigr\} \exp\bigl(-C_K(n\epsilon )^{1/3}\bigr)
\\
&=&  1-2C_0 \max\biggl\{\biggl(\frac{c_1p_1}{m_1}
\biggr)^{c_1m_1}, \biggl(\frac
{4c_1p_1}{m_1}\biggr)^{c_1m_1/4}\biggr\}
\exp\bigl(-C_K(n\epsilon)^{1/3}\bigr).
\end{eqnarray*}
Therefore, as $n>\epsilon^{-1}\max \{C_1, (c_1m_1\tilde
{C}_K^{-1}\log(c_1m_1 p_1))^3,(4^{-1}\tilde{C}_K^{-1}c_1m_1\log
(4c_1m_1 p_1))^3 \}$ where $\tilde{C}_K=C_K-\delta_0$ for some
$\delta_0>0$ being small, we have
%
%eA.15 #&#
\begin{equation}
\label{tilp0} P(\mathcal{A}_1\cap\mathcal{A}_2)>1-2C_0
\exp\bigl(-n^{1/3}\epsilon ^{1/3}\delta _0\bigr):=
\tilde{p}_{0n}\rightarrow1.
\end{equation}
\end{longlist}

%sA #&#
\section{Some auxiliary results}\label{appB}
%sA.1 #&#
\subsection{On the coverage of $\mathbb{B}_1$}
In this section, we consider the coverage of the set $\mathbb{B}_1$
used in the proof of Proposition~\ref{RENnorm}. For simplicity, for
any finite set $A$, let $\operatorname{conv} A$ stand for the convex hull of
$A$ and for any constant $C$, $C \operatorname{conv}(A)=\{C\cdot x \dvt x\in
\operatorname{conv}A \}$.
For $q = 1, 2$ and any positive integer $m$,
define $B^m_q = \{x \dvt \|x\|_q \leq 1, x \in R^m \}$.

%lethm.1 #&#
\begin{lemma}[(Mendelson, Pajor and  Tomczak-Jaegermann \cite{Mendelsonetal2008})]\label{Mendel1}
Let $m\geq1$ and $\varepsilon>0$. There exists an $\varepsilon$
cover $\Lambda\subset B_2^m$ of $B_2^m$ with respect to the Euclidean
metric such that
$B_2^m\subset(1-\varepsilon)^{-1} \operatorname{conv} \Lambda$ and
$|\Lambda|\leq(1+2/\varepsilon)^m$. Similarly there exists $\Lambda
'\subset S^{p-1}$ which is an $\varepsilon$ cover of the sphere
$S^{m-1}$ and $|\Lambda'|\leq(1+2/\varepsilon)^m$.
\end{lemma}

Define $\tilde{U}_m=\{x \in B_2^p \dvt |\operatorname{supp} x|\le m\}$ and
$U_m=\{x\in S^{p-1} \dvt |\operatorname{supp} x|\leq m\}$. The following result holds.

%lethm.2 #&#
\begin{lemma}[(Mendelson, Pajor and  Tomczak-Jaegermann \cite{Mendelsonetal2008})]\label{Mendel2}
There exists an
absolute constant $c$ for which the following holds.\vspace*{1pt} For every
$0<\varepsilon\leq1/2$ and every $1\leq m\leq p$,
there is a set $\Lambda\subset B_2^p$ which is an $\varepsilon$ cover
of $\tilde{U}_m$, such that $\tilde{U}_m\subset2 \operatorname{conv} \Lambda$ and $|\Lambda|$
is at most
%
%eA.1 #&#
\begin{equation}
\label{medcn}
\exp \biggl(cm\log \biggl(\frac{cp}{m\varepsilon} \biggr) \biggr).
\end{equation}
Moreover, there exists an $\varepsilon$ cover $\Lambda'\subset
S^{p-1}$ of $U_m$ with cardinality at most \textup{(\ref{medcn})}
Furthermore, for any $0\leq r\leq1$ there exists $\bar{\Lambda
}\subset r B_2^p$ such that $(U_m-U_m)\cap rB_2^p \subset2 \operatorname{conv} \bar{\Lambda}$ and $|\bar{\Lambda}|$ is at most \textup{(\ref{medcn})}.
\end{lemma}

For any $x\in R^{p_1}$, let $T_{0,x}$ with $|T_{0,x}|=s$ be the index
of the largest $s$ elements of $|x|=(|x_1|,\ldots,|x_{p_1}|)^T$.
Based on the results of Mendelson, Pajor and  Tomczak-Jaegermann \cite{Mendelsonetal2008}, we consider the
$\epsilon$-cover of the set
\begin{eqnarray}
\mathbb{B}_1&=&\bigl\{ x\in R^{p_1} \dvt \exists T
\mbox{ with } |T|\le s, \mbox{ such that } \|x_{T^c}\|_1
\leq k_0 \|x_{T}\|_1\bigr\}\cap
S^{p_1-1}
\nonumber
\\
&=&\bigl\{ x\in R^{p_1} \dvt \exists T \mbox{ with } |T|\le s, \mbox{ such
that } \|x_{T^c}\|_1 \leq k_0
\|x_{T}\|_1\bigr\}\cap\bigl\{x \dvt \|x\|_2\le1 \bigr\}\cap
S^{p_1-1}
\nonumber
\\
&:=&A_0\cap\bigl\{x\dvt \|x\|_2\le1\bigr\}\cap S^{p_1-1}.
\nonumber
\end{eqnarray}
By the result of Rudelson and Zhou \cite{Zhou2011}, it is easy to see
\begin{eqnarray*}
A_0&=&\bigl\{x\dvt x\in R^{p_1}, \exists T \mbox{ with } |T|
\le s, \mbox{ such that } \|x_{T^c}\|_1 \leq
k_0 \|x_{T}\|_1 \bigr\}
\nonumber
\\
&=&\bigl\{x\dvt x\in R^{p_1}, \|x_{T_{0,x}^c}\|_1 \leq
k_0 \|x_{T_{0,x}}\|_1 \bigr\}.
\end{eqnarray*}
Therefore, we consider the
coverage of
\begin{eqnarray*}
\mathbb{B}_1&=&\bigl\{ x\in R^{p_1} \dvt \|x_{T_{0,x}^c}
\|_1 \leq k_0 \| x_{T_{x,0}}\|_1, \|x
\|_2\le1 \bigr\}\cap S^{p_1-1}
\nonumber
\\
&:=&\tilde{A}_0\cap S^{p_1-1}.
\nonumber
\end{eqnarray*}
For $1\leq m\leq p_1$, let
\[
U_{m}=\bigl\{u\in S^{p_1-1} \dvt \operatorname{supp}(u)\leq m\bigr\}.
\]

%lethm.3 #&#
\begin{lemma}\label{lemm5}
The following two conclusions hold. (1) There exists $U_{m_1}$ of
$\mathbb{B}_1$, such that $\mathbb{B}_1 \subseteq2\operatorname{conv}
(U_{m_1})$ with $m_1=16s[(1+k_0)]^2$;
in addition there exists $\Lambda_1\subset S^{p_1-1}$ such that
$U_{m_1}\subset2 \operatorname{conv} \Lambda_1$ and the cardinality of
$\Lambda_1$ is given in \textup{(\ref{cover1})}.
(2) There exists $\Lambda_2\subseteq\epsilon B_2^{p_1}$, such that
$(\mathbb{B}_1-\mathbb{B}_1)\cap\epsilon B_{2}^{p_1}\subseteq
4\operatorname{conv}\Lambda_2$, for any $0\le\epsilon\le1$. And the
cardinality of $\Lambda_2$ is given in \textup{(\ref{cover2})}.
\end{lemma}

\begin{pf} \textit{Step} 1.  The proof of conclusion (1).
For $x\in\tilde{A}_0$, we have
\[
\|x\|_1\leq(1+k_0) \|x_{T_0}
\|_1 \leq s^{1/2}(1+k_0)\|x_{T_0}
\|_2\le s^{1/2}(1+k_0):=q(s).
\]
That is, $x\in q(s) B_1^{p_1}$.
Therefore, we have
\[
\mathbb{B}_1\subset q(s) B_1^{p_1}\cap
S^{p_1-1}:=\tilde{\mathbb{B}}_1.
\]
Therefore by Lemma~3.8 of Mendelson, Pajor and  Tomczak-Jaegermann \cite{Mendelsonetal2008},
it follows that $U_{m_1}$ with $m_1=4s[(1+k_0)/\epsilon]^2$ is a
$\epsilon$-cover of $\tilde{\mathbb{B}}_1$, for any $0<\epsilon
\le1/2$.
Taking $\epsilon=1/2$ leads to $m_1=16s[(1+k_0)]^2$.
By Lemma~\ref{Mendel2}, there is a set $\Lambda_1\subset S^{p_1-1}$,
such that $U_{m_1}\subset2\operatorname{conv} \Lambda_1 $ and
%
%eA.2 #&#
\begin{equation}
\label{cover1} |\Lambda_1|\leq\exp \bigl(c_1
m_1 \log [c_1p_1/m_1 ] \bigr)
\end{equation}
for some constant $c_1>0$. Therefore, we have $\mathbb{B}_1\subset
4\operatorname{conv} \Lambda_1$.

\textit{Step} 2. The proof of the second conclusion.

Since
$\mathbb{B}_1\subset\tilde{\mathbb{B}}_1$ and $\tilde{\mathbb{B}}_1$ is convex and star shaped, for any $0< \epsilon\le
1$, it holds that
\[
(\mathbb{B}_1-\mathbb{B}_1)\cap\epsilon
B_{2}^{p_1} \subseteq (\tilde{\mathbb{B}}_1-
\tilde{\mathbb{B}}_1)\cap\epsilon B_{2}^{p_1}
\subseteq \bigl[2q(s)B_1^{p_1}\cap2B_2^{p_1}
\bigr]\cap\epsilon B_{2}^{p_1}\subseteq 2q(s)B_1^{p_1}
\cap B_2^{p_1}.
\]
By Lemma~3.7 of Mendelson, Pajor and  Tomczak-Jaegermann \cite{Mendelsonetal2008}, for $m_2=4q^2(s)=4s(1+k_0)^2$,
we have
\[
(\mathbb{B}_1-\mathbb{B}_1)\cap
B_{2}^{p_1}\subseteq2 \operatorname{conv}(\tilde{U}_{m_2}).
\]
By Lemma~\ref{Mendel2}, there is a set $\Lambda_2'\subset B_2^{p_1}$
such that $\tilde{U}_{m_2}\subset2\operatorname{conv} \Lambda_2'$ and
\[
\bigl|\Lambda_2'\bigr|\leq\exp\bigl(c_1m_2
\log(c_1p_1/m_2)\bigr).
\]
Therefore, letting $\Lambda_2=\epsilon\Lambda_2'$ yields $(\mathbb
{B}_1-\mathbb{B}_1)\cap\epsilon B_{2}^{p_1}\subset4\operatorname{conv}
\Lambda_2$ and
$|\Lambda_2|=|\Lambda_2'|\leq\exp(c_1m_2\log[c_1p_1/m_2])$.
Recalling the definition of $m_1$ and $m_2$, we have $m_2=m_1/4$ and
%
%eA.3 #&#
\begin{equation}
\label{cover2} |\Lambda_2|=\bigl|\Lambda_2'\bigr|\leq
\exp\bigl(4^{-1}c_1m_1 \log[4c_1p_1/m_1]
\bigr).
\end{equation}
\end{pf}

%sA.2 #&#
\subsection{Concentration inequality of the squares of subexponential
type~variables}

%lethm.4 #&#
\begin{lemma}\label{lemm6}
Let $Z\in R$ be a subexponential type variable that satisfies
$P(|Z|>t)<C\exp(-t/K)$ as $t> C_0$ for some positive constants $C$, $K$ and $C_0$.
Assume that $E(Z)=0$ and that both $\operatorname{var}(Z)$ and
$\operatorname{var}(Z^2)$ are
finite. Let $Z_i, i=1,\ldots, n$ be the i.i.d. realizations of $Z$.
Then, for any $\delta>0$, as $n\delta>C_1$ for some $C_1>0$, we have
\[
P \Biggl( \Biggl|\frac{1}{n}\sum_{i=1}^n
Z_i^2-\operatorname {var}(Z) \Biggr|>\delta \Biggr)\le
C_2 \exp\bigl(-C_K(n\delta)^{1/3}\bigr),
\]
where $C_2>0$ and $C_K$ is a constant depending only on $K$.
\end{lemma}

%rethm.1 #&#
\begin{rem}
Similar to the subexponential distribution, we still denote $\|Z\|
_{\psi_1}=K$. From Lemma~\ref{lemm6}, we see that, for any fixed
$\delta>0$, as $n>C_1/\delta$, the inequality in Lemma~\ref{lemm6}
holds. Equivalently, for fixed $n$, if $\delta>C_1/n$, the above
inequality holds.
\end{rem}

\begin{pf*}{Proof of Lemma~\protect\ref{lemm6}}
Without loss of generality, we assume $\operatorname{var}(Z)=1$
in the
following argument. Let $A_n=n^{-1}\sum_{i=1}^nZ^2_i-1$. For
any fixed large constant $M>\max(C_0,1)$, let $Y_i^M=(Z_i^2-1)I_{\{
|Z_i|<M\}}$ and $X_i^M=(Z_i^2-1)I_{\{|Z_i|>M\}}$. Then we have
\[
A_n=n^{-1}\sum_{i=1}^n
Y_i^M+n^{-1}\sum_{i=1}^n
X_i^M:=A_{1n}+A_{2n}.
\]
Note that $E(A_n)=E(A_{1n})+E(A_{2n})=0$. It holds that
\[
A_n=A_{1n}-E(A_{1n})+A_{2n}-E(A_{2n}).
\]
By the Bernstein inequality, due to the fact that $|Y_i^M|\le M^2-1$,
we have, for any $\delta>0$
%
%eA.4 #&#
\begin{eqnarray}
P\bigl(\bigl|A_{1n}-E(A_{1n})\bigr|>\delta/2\bigr)&=&P
\Biggl(\Biggl|n^{-1}\sum_{i=1}^n
Y_i^M-E\bigl(Y_i^M\bigr)\Biggr|>
\delta/2\Biggr)
\nonumber
\\
&\le&2\exp \biggl(-\frac{(n\delta/2)^2/2}{nE(Y_1^M)^2+(M^2-1)n\delta
/6} \biggr)
\nonumber
\\[-8pt]
\label{P1}
\\[-8pt]
\nonumber
&\le&2 \exp \biggl(-\frac{n\delta^2}{8E(Y_1^M)^2+4(M^2-1)\delta/3} \biggr)
\nonumber
\\
&\le&2\exp \biggl(-\frac{n\delta^2}{8c_0+4(M^2-1)\delta/3} \biggr),
\nonumber
\end{eqnarray}
where $c_0=\operatorname{var} (Z^2) <\infty$.
In addition,
\[
P\bigl(\bigl|A_{2n}-E(A_{2n})\bigr|>\delta/2\bigr)\le
\frac{4}{n\delta^2}\operatorname {var}\bigl(X_1^M\bigr)\le
\frac{4}{n\delta^2}E\bigl(Z_1^4I_{\{|Z_1|>M\}}\bigr).
\]
Let $V=|Z_1|$ and note that $V$ is subexponential. By integration by
parts, it is easy to see
\begin{eqnarray}
E\bigl(V^4I_{\{V>M\}}\bigr)&=&P(V>M)M^4+ \int
_{x>M} 4P(V>x) x^3\,\mathrm{d}x.
\nonumber
\end{eqnarray}
Because $V$ is subexponential, it follows that $P(V>s)\le C\exp(-
s/K)$ for any $s>M$. Inserting this inequality into the above equation
and applying integration by part to the second term repeatedly, for
large $M$, we have
\[
E\bigl(V^4I_{\{V>M\}}\bigr)\le c_1M^4
\exp(-C_K M),
\]
for some  $c_1>0$
where $C_K$ is a constant depending only on $K$.
Therefore, we have
%
%eA.5 #&#
\begin{eqnarray}
\label{P2} P\bigl(\bigl|A_{2n}-E(A_{2n})\bigr|>\delta/2\bigr)\le
\frac{c_1M^4}{n\delta^2}\exp(-C_K M).
\end{eqnarray}
For $M=c_{2}(n\delta)^{1/3}$ for some constant $c_2$, the right-hand
side of (\ref{P1}) has the same order as that of (\ref{P2}).
Also, we remind that we need $M>\max(C_0,1)$. Therefore, we have
$n\delta>[c_2^{-1}\max(C_0,1)]^{3}:=C_1$, that is $\delta$ cannot
be too small. Consequently, for any $\delta>C_1/n$, we have
\[
P\bigl(|A_n|>\delta\bigr)\le P\bigl(\bigl|A_{1n}-E(A_{1n})\bigr|>
\delta /2\bigr)+P\bigl(\bigl|A_{2n}-E(A_{2n})\bigr|>\delta/2\bigr)\le
C_2 \exp\bigl[-C_K (n\delta)^{1/3}\bigr]
\]
for some constant $C_2>0$ and $C_K$, where $C_K$ depends only on $K$.
The proof is completed.
\end{pf*}

%sA #&#
\section{Some properties of sub-Gaussian and
subexponential distributions}\label{appC}

The following two propositions on sub-Gaussian and subexponential
distributions are from Vershynin \cite{Vershynin2011}. They are listed here for
completeness.

%prthm.1 #&#
\begin{prop}[(Bernstein-type inequality)]\label{ApC1} Let $X_1, \ldots
,X_n$ be independent centered subexponential random variables, and $K =
\max_i \|X_i\|_{\psi_1}$. Then for every $a =(a_1, \ldots, a_n)\in
R^n$ and every $t \ge0$, we have
\[
P\Biggl(\Biggl|\sum_{i=1}^n a_iX_i
\Biggr|\ge t \Biggr)\le2 \exp \biggl(- c\min\biggl\{ \frac{t^2}{K^2\|a\|_2^2},
\frac{t}{K \|a\|_\infty}\biggr\} \biggr),
\]
where $c> 0$ is an absolute constant.
\end{prop}

%prthm.2 #&#
\begin{prop}\label{ApC2}  A random variable $X$ is sub-Gaussian if
and only if $X^2$ is sub-exponential. Moreover,
$\|X\|^2_{\psi_2} \le\|X^2\|_{\psi_1} \le2\|X\|^2_{\psi_2}$.
\end{prop}

\begin{cor}
If $X, Y$ are centered
sub-Gaussian with $\max\{\|X\|_{\psi_2}, \|Y\|_{\psi_2}\}<K$. Then
$XY$ is subexponential with $\|XY\|_{\psi_1}\le2K^2$.
\end{cor}

\begin{pf}
For any $t>0$, $P(|XY|>t)\le P(X^2+Y^2>2t)\le
P(X^2>t)+P(Y^2>t)$. By\vspace*{1pt} Proposition~\ref{ApC2}, we have that $X^2$ and
$Y^2$ are both subexponential with $\psi_{1}$ norm less than $2K^2$;
that is, $P(X^2>t)\le c \exp(-t/2K^2)$ for some constant $c>0$ and the
same is true for $Y^2$.
Consequently, we have $P(|XY|>t)\le2c\exp(-t/2K^2)$. The conclusion holds.
\end{pf}
\end{appendix}

\section*{Acknowledgements}

We are grateful to the Associate Editor and an anonymous referee
for their helpful comments. The researches of Zhao is supported by
National Science Foundation of China (Nos. 11101022, 11471030) and
Fundamental Research Funds for the Central Universities.

% imsref loaded by daiva.urboniene, 2015-05-21 09:55:59

\printhistory

\begin{thebibliography}{33}
% pybtex-1.32. Style name=bej, version=1.42, label_style=nameyear, sorting_style=complex, cfg=None, language=None.


%b1 ###
%b1 #&#
\bibitem{Bachetal2012}
\begin{barticle}[mr]
\bauthor{\bsnm{Bach},~\bfnm{Francis}\binits{F.}},
\bauthor{\bsnm{Jenatton},~\bfnm{Rodolphe}\binits{R.}},
\bauthor{\bsnm{Mairal},~\bfnm{Julien}\binits{J.}} \AND
\bauthor{\bsnm{Obozinski},~\bfnm{Guillaume}\binits{G.}}
(\byear{2012}).
\btitle{Structured sparsity through convex optimization}.
\bjournal{Statist. Sci.}
\bvolume{27}
\bpages{450--468}.
\bid{doi={10.1214/12-STS394}, issn={0883-4237}, mr={3025128}}
\end{barticle}
%
\iffalse\OrigBibText
Bach, F., Jenatton, R., Mairal, J. and
Obozinski, G. (2012).
Structured sparsity through convex optimization. \textit{Stat. Sci.} \textbf{27}, 450--468.
\endOrigBibText\fi
\bptok{imsref}%
% NOT OUTPUTTED:
%   number = 4
%   doi = http://dx.doi.org/10.1214/12-STS394
%   fjournal = Statistical Science. A Review Journal of the Institute of Mathematical Statistics
\endbibitem

%b2 ###
%b2 #&#
\bibitem{Bickeletal2009}
\begin{barticle}[mr]
\bauthor{\bsnm{Bickel},~\bfnm{Peter~J.}\binits{P.J.}},
\bauthor{\bsnm{Ritov},~\bfnm{Ya'acov}\binits{Y.}} \AND
\bauthor{\bsnm{Tsybakov},~\bfnm{Alexandre~B.}\binits{A.B.}}
(\byear{2009}).
\btitle{Simultaneous analysis of lasso and {D}antzig selector}.
\bjournal{Ann. Statist.}
\bvolume{37}
\bpages{1705--1732}.
\bid{doi={10.1214/08-AOS620}, issn={0090-5364}, mr={2533469}}
\end{barticle}
%
\iffalse\OrigBibText
Bickel, P., Ritov, Y. and Tsybakov, A.
(2009).
Simultaneous analysis of Lasso and Dantzig selector.
\textit{Ann. Statist.} \textbf{37}, 1705--1732.
\endOrigBibText\fi
\bptok{imsref}%
% NOT OUTPUTTED:
%   number = 4
%   doi = http://dx.doi.org/10.1214/08-AOS620
%   coden = ASTSC7
%   fjournal = The Annals of Statistics
\endbibitem

%b3 ###
%b3 #&#
\bibitem{Bickeletal2010}
\begin{bincollection}[mr]
\bauthor{\bsnm{Bickel},~\bfnm{Peter~J.}\binits{P.J.}},
\bauthor{\bsnm{Ritov},~\bfnm{Ya'acov}\binits{Y.}} \AND
\bauthor{\bsnm{Tsybakov},~\bfnm{Alexandre~B.}\binits{A.B.}}
(\byear{2010}).
\btitle{Hierarchical selection of variables in sparse high-dimensional regression}.
In \bbooktitle{Borrowing Strength: Theory Powering Applications---a {F}estschrift for {L}awrence D. {B}rown}.
\bseries{Inst. Math. Stat. Collect.}
\bvolume{6}
\bpages{56--69}.
\blocation{Beachwood, OH}:
\bpublisher{IMS}.
\bid{mr={2798511}}
\end{bincollection}
%
\iffalse\OrigBibText
Bickel, P., Ritov, Y. and Tsybakov, A. (2010).
Hierarchical selection of variables in sparse high-dimensional
regression. \textit{IMS Collections.
Borrowing Strength: Theory Powering Applications: Festschrift for
Lawrence D. Brown} \textbf{6}, 56--69.
\endOrigBibText\fi
\bptok{imsref}%
\endbibitem

%b4 ###
%b4 #&#
\bibitem{Bienetal2013}
\begin{barticle}[mr]
\bauthor{\bsnm{Bien},~\bfnm{Jacob}\binits{J.}},
\bauthor{\bsnm{Taylor},~\bfnm{Jonathan}\binits{J.}} \AND
\bauthor{\bsnm{Tibshirani},~\bfnm{Robert}\binits{R.}}
(\byear{2013}).
\btitle{A LASSO for hierarchical interactions}.
\bjournal{Ann. Statist.}
\bvolume{41}
\bpages{1111--1141}.
\bid{doi={10.1214/13-AOS1096}, issn={0090-5364}, mr={3113805}}
\end{barticle}
%
\iffalse\OrigBibText
Bien, J., Taylor, J. and Tibshirani, R.
(2013). A Lasso for hierarchical interactions.
\textit{Ann. Statist.} \textbf{41}, 1111--1141.
\endOrigBibText\fi
\bptok{imsref}%
% NOT OUTPUTTED:
%   number = 3
%   doi = http://dx.doi.org/10.1214/13-AOS1096
%   fjournal = The Annals of Statistics
\endbibitem

%b5 ###
%b5 #&#
\bibitem{CandesTao2007}
\begin{barticle}[mr]
\bauthor{\bsnm{Candes},~\bfnm{Emmanuel}\binits{E.}} \AND
\bauthor{\bsnm{Tao},~\bfnm{Terence}\binits{T.}}
(\byear{2007}).
\btitle{The {D}antzig selector: Statistical estimation when {$p$} is much larger than~{$n$}}.
\bjournal{Ann. Statist.}
\bvolume{35}
\bpages{2313--2351}.
\bid{doi={10.1214/009053606000001523}, issn={0090-5364}, mr={2382644}}
\bptnote{check related}%
\end{barticle}
%
\iffalse\OrigBibText
Candes, E. and Tao, T. (2007). The Dantzig selector: Statistical
estimation when $p$ is much larger than $n$ (with discussion).
\textit{Ann. Statist.} \textbf{35}, 2313--2351.
\endOrigBibText\fi
\bptok{imsref}%
% NOT OUTPUTTED:
%   number = 6
%   doi = http://dx.doi.org/10.1214/009053606000001523
%   coden = ASTSC7
%   fjournal = The Annals of Statistics
\endbibitem

%b6 ###
%b6 #&#
\bibitem{Choietal2010}
\begin{barticle}[mr]
\bauthor{\bsnm{Choi},~\bfnm{Nam~Hee}\binits{N.H.}},
\bauthor{\bsnm{Li},~\bfnm{William}\binits{W.}} \AND
\bauthor{\bsnm{Zhu},~\bfnm{Ji}\binits{J.}}
(\byear{2010}).
\btitle{Variable selection with the strong heredity constraint and its oracle property}.
\bjournal{J. Amer. Statist. Assoc.}
\bvolume{105}
\bpages{354--364}.
%\bnote{With supplementary material available online}.
\bid{doi={10.1198/jasa.2010.tm08281}, issn={0162-1459}, mr={2656056}}
\end{barticle}
%
\iffalse\OrigBibText
Choi, N. H., Li, W. and Zhu, J. (2010). Variable selection with the strong
heredity constraint and its oracle property. \textit{ J. Am. Statist.
Assoc.} \textbf{105}, 354--364.
\endOrigBibText\fi
\bptok{imsref}%
% NOT OUTPUTTED:
%   number = 489
%   doi = http://dx.doi.org/10.1198/jasa.2010.tm08281
%   coden = JSTNAL
%   fjournal = Journal of the American Statistical Association
\endbibitem

%b7 ###
%b7 #&#
\bibitem{FanLi2001}
\begin{barticle}[mr]
\bauthor{\bsnm{Fan},~\bfnm{Jianqing}\binits{J.}} \AND
\bauthor{\bsnm{Li},~\bfnm{Runze}\binits{R.}}
(\byear{2001}).
\btitle{Variable selection via nonconcave penalized likelihood and its oracle properties}.
\bjournal{J. Amer. Statist. Assoc.}
\bvolume{96}
\bpages{1348--1360}.
\bid{doi={10.1198/016214501753382273}, issn={0162-1459}, mr={1946581}}
\end{barticle}
%
\iffalse\OrigBibText
Fan, J. and Li, R. (2001).
Variable selection via nonconcave penalized likelihood and its
oracle properties.
\textit{J. Am. Statist. Assoc.} \textbf{96}, 1348--1360.
\endOrigBibText\fi
\bptok{imsref}%
% NOT OUTPUTTED:
%   number = 456
%   doi = http://dx.doi.org/10.1198/016214501753382273
%   coden = JSTNAL
%   fjournal = Journal of the American Statistical Association
\endbibitem

%b8 ###
%b8 #&#
\bibitem{Fangetal1994}
\begin{barticle}[mr]
\bauthor{\bsnm{Fang},~\bfnm{Yu~Guang}\binits{Y.G.}},
\bauthor{\bsnm{Loparo},~\bfnm{Kenneth~A.}\binits{K.A.}} \AND
\bauthor{\bsnm{Feng},~\bfnm{Xiangbo}\binits{X.}}
(\byear{1994}).
\btitle{Inequalities for the trace of matrix product}.
\bjournal{IEEE Trans. Automat. Control}
\bvolume{39}
\bpages{2489--2490}.
\bid{doi={10.1109/9.362841}, issn={0018-9286}, mr={1337578}}
\end{barticle}
%
\iffalse\OrigBibText
Fang, Y., Loparo, K. A. and Feng, X. (1994).
Inequalities for
the trace of matrix product.
\textit{IEEE Trans. Autom. Control},
\textbf{39}, 2489--2490.
\endOrigBibText\fi
\bptok{imsref}%
% NOT OUTPUTTED:
%   number = 12
%   doi = http://dx.doi.org/10.1109/9.362841
%   coden = IETAA9
%   fjournal = Institute of Electrical and Electronics Engineers. Transactions on Automatic Control
\endbibitem

%b9 ###
%b9 #&#
\bibitem{HallXue2013}
\begin{barticle}[mr]
\bauthor{\bsnm{Hall},~\bfnm{Peter}\binits{P.}} \AND
\bauthor{\bsnm{Xue},~\bfnm{Jing-Hao}\binits{J.-H.}}
(\byear{2014}).
\btitle{On selecting interacting features from high-dimensional data}.
\bjournal{Comput. Statist. Data Anal.}
\bvolume{71}
\bpages{694--708}.
\bid{doi={10.1016/j.csda.2012.10.010}, issn={0167-9473}, mr={3132000}}
\end{barticle}
%
\iffalse\OrigBibText
Hall, P. and Xue, J. H. (2014). On selecting interacting features from
high-dimensional data.
\textit{Comput. Stat. Data An.} \textbf{71}, 694--708.
\endOrigBibText\fi
\bptok{imsref}%
% NOT OUTPUTTED:
%   doi = http://dx.doi.org/10.1016/j.csda.2012.10.010
%   fjournal = Computational Statistics \& Data Analysis
\endbibitem

%b10 ###
%b10 #&#
\bibitem{HaoZhang2012a}
\begin{bmisc}[auto:parserefs-M02]
\bauthor{\bsnm{Hao},~\bfnm{N.}\binits{N.}} \AND
\bauthor{\bsnm{Zhang},~\bfnm{H.~H.}\binits{H.H.}}
(\byear{2012}).
\bhowpublished{Interaction selection under marginality principle in high dimensional regression.
Manuscript.}
\end{bmisc}
%
\iffalse\OrigBibText
Hao, N. and Zhang, H. H. (2012a). Interaction selection under
marginality principle in high dimensional regression. Manuscript.
\endOrigBibText\fi
\bptok{imsref}%
% NOT OUTPUTTED:
%   sortkey = Hao(2012
\endbibitem

%b11 ###
%b11 #&#
\bibitem{HaoZhang2012b}
\begin{bmisc}[auto:parserefs-M02]
\bauthor{\bsnm{Hao},~\bfnm{N.}\binits{N.}} \AND
\bauthor{\bsnm{Zhang},~\bfnm{H.~H.}\binits{H.H.}}
(\byear{2012}).
\bhowpublished{A note on regression models with interactions.
Manuscript.}
\end{bmisc}
%
\iffalse\OrigBibText
Hao, N. and Zhang, H. H. (2012b). A note on regression models with
interactions. Manuscript.
\endOrigBibText\fi
\bptok{imsref}%
% NOT OUTPUTTED:
%   sortkey = Hao(2012
\endbibitem

%b12 ###
%b12 #&#
\bibitem{HaoZhang2014}
\begin{barticle}[mr]
\bauthor{\bsnm{Hao},~\bfnm{Ning}\binits{N.}} \AND
\bauthor{\bsnm{Zhang},~\bfnm{Hao~Helen}\binits{H.H.}}
(\byear{2014}).
\btitle{Interaction screening for ultrahigh-dimensional data}.
\bjournal{J. Amer. Statist. Assoc.}
\bvolume{109}
\bpages{1285--1301}.
\bid{doi={10.1080/01621459.2014.881741}, issn={0162-1459}, mr={3265697}}
\end{barticle}
%
\iffalse\OrigBibText
Hao, N. and Zhang, H. H. (2014). Interaction screening for ultra-high
dimensional data. \textit{J. Am. Statist. Assoc.} \textbf{109}, 1285--1301.
\endOrigBibText\fi
\bptok{imsref}%
% NOT OUTPUTTED:
%   number = 507
%   doi = http://dx.doi.org/10.1080/01621459.2014.881741
%   fjournal = Journal of the American Statistical Association
\endbibitem

%b13 ###
%b13 #&#
\bibitem{LinZhang2006}
\begin{barticle}[mr]
\bauthor{\bsnm{Lin},~\bfnm{Yi}\binits{Y.}} \AND
\bauthor{\bsnm{Zhang},~\bfnm{Hao~Helen}\binits{H.H.}}
(\byear{2006}).
\btitle{Component selection and smoothing in multivariate nonparametric regression}.
\bjournal{Ann. Statist.}
\bvolume{34}
\bpages{2272--2297}.
\bid{doi={10.1214/009053606000000722}, issn={0090-5364}, mr={2291500}}
\end{barticle}
%
\iffalse\OrigBibText
Lin, Y. and Zhang, H. H. (2006). Component selection and smoothing in
multivariate non-parametric Regression.
\textit{Ann. Statist.}, \textbf{34}, 2272--2297.
\endOrigBibText\fi
\bptok{imsref}%
% NOT OUTPUTTED:
%   number = 5
%   doi = http://dx.doi.org/10.1214/009053606000000722
%   coden = ASTSC7
%   fjournal = The Annals of Statistics
\endbibitem

%b14 ###
%b14 #&#
\bibitem{MagnusNeudecker1979}
\begin{barticle}[mr]
\bauthor{\bsnm{Magnus},~\bfnm{Jan~R.}\binits{J.R.}} \AND
\bauthor{\bsnm{Neudecker},~\bfnm{H.}\binits{H.}}
(\byear{1979}).
\btitle{The commutation matrix: Some properties and applications}.
\bjournal{Ann. Statist.}
\bvolume{7}
\bpages{381--394}.
\bid{issn={0090-5364}, mr={0520247}}
\end{barticle}
%
\iffalse\OrigBibText
Magnus, J. and Neudecker, H.
(1979). The commutation matrix: Some properties and applications.
\textit{Ann. Statist.} \textbf{7}, 381--394.
\endOrigBibText\fi
\bptok{imsref}%
% NOT OUTPUTTED:
%   url = http://links.jstor.org/sici?sici=0090-5364(197903)7:2<381:TCMSPA>2.0.CO;2-Y&origin=MSN
%   number = 2
%   coden = ASTSC7
%   fjournal = The Annals of Statistics
\endbibitem

%b15 ###
%b15 #&#
\bibitem{MeinshausenYu2009}
\begin{barticle}[mr]
\bauthor{\bsnm{Meinshausen},~\bfnm{Nicolai}\binits{N.}} \AND
\bauthor{\bsnm{Yu},~\bfnm{Bin}\binits{B.}}
(\byear{2009}).
\btitle{Lasso-type recovery of sparse representations for high-dimensional data}.
\bjournal{Ann. Statist.}
\bvolume{37}
\bpages{246--270}.
\bid{doi={10.1214/07-AOS582}, issn={0090-5364}, mr={2488351}}
\end{barticle}
%
\iffalse\OrigBibText
Meinshausen, N. and Yu, B. (2009). Lasso type recovery of sparse
representations for high dimensional data. \textit{Ann. Statist.}
\textbf{37}, 246--270.
\endOrigBibText\fi
\bptok{imsref}%
% NOT OUTPUTTED:
%   number = 1
%   doi = http://dx.doi.org/10.1214/07-AOS582
%   coden = ASTSC7
%   fjournal = The Annals of Statistics
\endbibitem

%b16 ###
%b16 #&#
\bibitem{Mendelsonetal2008}
\begin{barticle}[mr]
\bauthor{\bsnm{Mendelson},~\bfnm{Shahar}\binits{S.}},
\bauthor{\bsnm{Pajor},~\bfnm{Alain}\binits{A.}} \AND
\bauthor{\bsnm{Tomczak-Jaegermann},~\bfnm{Nicole}\binits{N.}}
(\byear{2008}).
\btitle{Uniform uncertainty principle for {B}ernoulli and subgaussian ensembles}.
\bjournal{Constr. Approx.}
\bvolume{28}
\bpages{277--289}.
\bid{doi={10.1007/s00365-007-9005-8}, issn={0176-4276}, mr={2453368}}
\end{barticle}
%
\iffalse\OrigBibText
Mendelson, S., Pajor, A., and
Tomczak-Jaegermann, N. (2008). {Uniform uncertainty principle for
Bernoulli and sub-Gaussian ensembles}.
\textit{Constr. Approx.} \textbf{28}, 277--289.
\endOrigBibText\fi
\bptok{imsref}%
% NOT OUTPUTTED:
%   number = 3
%   doi = http://dx.doi.org/10.1007/s00365-007-9005-8
%   fjournal = Constructive Approximation. An International Journal for Approximations and Expansions
\endbibitem

%b17 ###
%b17 #&#
\bibitem{Negahban2012}
\begin{barticle}[mr]
\bauthor{\bsnm{Negahban},~\bfnm{Sahand~N.}\binits{S.N.}},
\bauthor{\bsnm{Ravikumar},~\bfnm{Pradeep}\binits{P.}},
\bauthor{\bsnm{Wainwright},~\bfnm{Martin~J.}\binits{M.J.}} \AND
\bauthor{\bsnm{Yu},~\bfnm{Bin}\binits{B.}}
(\byear{2012}).
\btitle{A unified framework for high-dimensional analysis of {$M$}-estimators with decomposable regularizers}.
\bjournal{Statist. Sci.}
\bvolume{27}
\bpages{538--557}.
\bid{doi={10.1214/12-STS400}, issn={0883-4237}, mr={3025133}}
\end{barticle}
%
\iffalse\OrigBibText
Negahban, S. N., Ravikumar, P., Wainwright,
M. J. and Yu, B. (2012). A Uuified framework for high-dimensional
analysis of $M$-estimators with decomposable regularizers. \textit{
Stat. Sci. } \textbf{27}, 538--557.
\endOrigBibText\fi
\bptok{imsref}%
% NOT OUTPUTTED:
%   number = 4
%   doi = http://dx.doi.org/10.1214/12-STS400
%   fjournal = Statistical Science. A Review Journal of the Institute of Mathematical Statistics
\endbibitem

%b18 ###
%b18 #&#
\bibitem{RadchenkoJames2010}
\begin{barticle}[mr]
\bauthor{\bsnm{Radchenko},~\bfnm{Peter}\binits{P.}} \AND
\bauthor{\bsnm{James},~\bfnm{Gareth~M.}\binits{G.M.}}
(\byear{2010}).
\btitle{Variable selection using adaptive nonlinear interaction structures in high dimensions}.
\bjournal{J. Amer. Statist. Assoc.}
\bvolume{105}
\bpages{1541--1553}.
\bid{doi={10.1198/jasa.2010.tm10130}, issn={0162-1459}, mr={2796570}}
\end{barticle}
%
\iffalse\OrigBibText
Radchenko, P. and James, G. M.
(2010). Variable selection using adaptive non-linear interaction
structures in high dimensions. \textit{J. Am. Statist. Assoc.}
\textbf{105}, 1541--1553.
\endOrigBibText\fi
\bptok{imsref}%
% NOT OUTPUTTED:
%   number = 492
%   doi = http://dx.doi.org/10.1198/jasa.2010.tm10130
%   coden = JSTNAL
%   fjournal = Journal of the American Statistical Association
\endbibitem

%b19 ###
%b19 #&#
\bibitem{Raskuttietal2010}
\begin{barticle}[mr]
\bauthor{\bsnm{Raskutti},~\bfnm{Garvesh}\binits{G.}},
\bauthor{\bsnm{Wainwright},~\bfnm{Martin~J.}\binits{M.J.}} \AND
\bauthor{\bsnm{Yu},~\bfnm{Bin}\binits{B.}}
(\byear{2010}).
\btitle{Restricted eigenvalue properties for correlated {G}aussian designs}.
\bjournal{J. Mach. Learn. Res.}
\bvolume{11}
\bpages{2241--2259}.
\bid{issn={1532-4435}, mr={2719855}}
\end{barticle}
%
\iffalse\OrigBibText
Raskutti, G., Wainwright, M. J. and Yu, B. (2010). Restricted
eigenvalue properties for correlated Gaussian designs. \textit{J.
Mach. Learn. Res.} \textbf{11}, 2241--2259.
\endOrigBibText\fi
\bptok{imsref}%
% NOT OUTPUTTED:
%   fjournal = Journal of Machine Learning Research (JMLR)
\endbibitem

%b20 ###
%b20 #&#
\bibitem{Raskutti2011}
\begin{barticle}[mr]
\bauthor{\bsnm{Raskutti},~\bfnm{Garvesh}\binits{G.}},
\bauthor{\bsnm{Wainwright},~\bfnm{Martin~J.}\binits{M.J.}} \AND
\bauthor{\bsnm{Yu},~\bfnm{Bin}\binits{B.}}
(\byear{2011}).
\btitle{Minimax rates of estimation for high-dimensional linear regression over {$\ell\sb q$}-balls}.
\bjournal{IEEE Trans. Inform. Theory}
\bvolume{57}
\bpages{6976--6994}.
\bid{doi={10.1109/TIT.2011.2165799}, issn={0018-9448}, mr={2882274}}
\end{barticle}
%
\iffalse\OrigBibText
Raskutti, G., Wainwright, M. J. and Yu, B.
(2011).
Minimax rates of estimation for high-dimensional linear
regression over-balls. \textit{IEEE Trans. Inf. Theory} \textbf{57},
6976-6994.
\endOrigBibText\fi
\bptok{imsref}%
% NOT OUTPUTTED:
%   number = 10
%   doi = http://dx.doi.org/10.1109/TIT.2011.2165799
%   coden = IETTAW
%   fjournal = Institute of Electrical and Electronics Engineers. Transactions on Information Theory
\endbibitem

%b21 ###
%b21 #&#
\bibitem{Zhou2011}
\begin{barticle}[mr]
\bauthor{\bsnm{Rudelson},~\bfnm{Mark}\binits{M.}} \AND
\bauthor{\bsnm{Zhou},~\bfnm{Shuheng}\binits{S.}}
(\byear{2013}).
\btitle{Reconstruction from anisotropic random measurements}.
\bjournal{IEEE Trans. Inform. Theory}
\bvolume{59}
\bpages{3434--3447}.
\bid{doi={10.1109/TIT.2013.2243201}, issn={0018-9448}, mr={3061256}}
\end{barticle}
%
\iffalse\OrigBibText
Rudelson, M. and Zhou, S. (2013).
Reconstruction from anisotropic random measurements. \textit{IEEE
Trans. Inf. Theory}, \textbf{59}, 3434-3447.
\endOrigBibText\fi
\bptok{imsref}%
% NOT OUTPUTTED:
%   number = 6
%   doi = http://dx.doi.org/10.1109/TIT.2013.2243201
%   coden = IETTAW
%   fjournal = Institute of Electrical and Electronics Engineers. Transactions on Information Theory
\endbibitem

%b22 ###
%b22 #&#
\bibitem{Shah2012}
\begin{bmisc}[auto:parserefs-M02]
\bauthor{\bsnm{Shah},~\bfnm{R.~D.}\binits{R.D.}}
(\byear{2012}).
\bhowpublished{Modelling interactions in high-dimensional data with backtracking.
Manuscript.}
\end{bmisc}
%
\iffalse\OrigBibText
Shah, R. D. (2012).
Modelling interactions in high-dimensional data
with backtracking. Manuscript.
\endOrigBibText\fi
\bptok{imsref}%
% NOT OUTPUTTED:
%   sortkey = Shah(2012
\endbibitem

%b23 ###
%b23 #&#
\bibitem{Tibshirani1996}
\begin{barticle}[mr]
\bauthor{\bsnm{Tibshirani},~\bfnm{Robert}\binits{R.}}
(\byear{1996}).
\btitle{Regression shrinkage and selection via the lasso}.
\bjournal{J. Roy. Statist. Soc. Ser. B}
\bvolume{58}
\bpages{267--288}.
\bid{issn={0035-9246}, mr={1379242}}
\end{barticle}
%
\iffalse\OrigBibText
Tibshirani, R. (1996). Regression shrinkage and selection via the
Lasso. \textit{J. Roy. Stat. Soc. B. Met.} \textbf{58}, 267--288.
\endOrigBibText\fi
\bptok{imsref}%
% NOT OUTPUTTED:
%   url = http://links.jstor.org/sici?sici=0035-9246(1996)58:1<267:RSASVT>2.0.CO;2-G&origin=MSN
%   number = 1
%   coden = JSTBAJ
%   fjournal = Journal of the Royal Statistical Society. Series B. Methodological
\endbibitem

%b24 ###
%b24 #&#
\bibitem{vandeGeer2007}
\begin{bmisc}[auto:parserefs-M02]
\bauthor{\bsnm{van~de Geer},~\bfnm{S.}\binits{S.}}
(\byear{2007}).
\bhowpublished{The deterministic lasso.
Techinical Report 140, ETH, Zurich.}
\end{bmisc}
%
\iffalse\OrigBibText
van de Geer, S. (2007). The deterministic lasso. Techinical
Report 140, ETH Zurich.
\endOrigBibText\fi
\bptok{imsref}%
\endbibitem

%b25 ###
%b25 #&#
\bibitem{GeerBuhl2009}
\begin{barticle}[mr]
\bauthor{\bsnm{van~de Geer},~\bfnm{Sara~A.}\binits{S.A.}} \AND
\bauthor{\bsnm{B{\"u}hlmann},~\bfnm{Peter}\binits{P.}}
(\byear{2009}).
\btitle{On the conditions used to prove oracle results for the {L}asso}.
\bjournal{Electron. J. Stat.}
\bvolume{3}
\bpages{1360--1392}.
\bid{doi={10.1214/09-EJS506}, issn={1935-7524}, mr={2576316}}
\end{barticle}
%
\iffalse\OrigBibText
van de Geer, S. and B\"uhlmann, P. (2009). On the conditions used to
prove oracle results for the Lasso. \textit{Electron. J. Stat.}
\textbf{3}, 1360--1392.
\endOrigBibText\fi
\bptok{imsref}%
% NOT OUTPUTTED:
%   doi = http://dx.doi.org/10.1214/09-EJS506
%   fjournal = Electronic Journal of Statistics
\endbibitem

%b26 ###
%b26 #&#
\bibitem{Vershynin2011}
\begin{bincollection}[mr]
\bauthor{\bsnm{Vershynin},~\bfnm{Roman}\binits{R.}}
(\byear{2012}).
\btitle{Introduction to the non-asymptotic analysis of random matrices}.
In \bbooktitle{Compressed Sensing}
\bpages{210--268}.
\blocation{Cambridge}:
\bpublisher{Cambridge Univ. Press}.
\bid{mr={2963170}}
\bptnote{check pages, check year}%
\end{bincollection}
%
\iffalse\OrigBibText
Vershynin, R. (2011).
Introduction to the non-asymptotic analysis of
random matrices. arXiv:1011.3027v5.
\endOrigBibText\fi
\bptok{imsref}%
\endbibitem

%b27 ###
%b27 #&#
\bibitem{Wang2009}
\begin{barticle}[mr]
\bauthor{\bsnm{Wang},~\bfnm{Hansheng}\binits{H.}}
(\byear{2009}).
\btitle{Forward regression for ultra-high dimensional variable screening}.
\bjournal{J. Amer. Statist. Assoc.}
\bvolume{104}
\bpages{1512--1524}.
\bid{doi={10.1198/jasa.2008.tm08516}, issn={0162-1459}, mr={2750576}}
\end{barticle}
%
\iffalse\OrigBibText
Wang, H. (2009).
Forward regression for ultra-high dimensional variable
screening. \textit{J. Am. Statist. Assoc.} \textbf{104}, 1512--1524.
\endOrigBibText\fi
\bptok{imsref}%
% NOT OUTPUTTED:
%   number = 488
%   doi = http://dx.doi.org/10.1198/jasa.2008.tm08516
%   coden = JSTNAL
%   fjournal = Journal of the American Statistical Association
\endbibitem

%b28 ###
%b28 #&#
\bibitem{YeZhang2010}
\begin{barticle}[mr]
\bauthor{\bsnm{Ye},~\bfnm{Fei}\binits{F.}} \AND
\bauthor{\bsnm{Zhang},~\bfnm{Cun-Hui}\binits{C.-H.}}
(\byear{2010}).
\btitle{Rate minimaxity of the {L}asso and {D}antzig selector for the {$\ell\sb q$} loss in {$\ell\sb r$} balls}.
\bjournal{J. Mach. Learn. Res.}
\bvolume{11}
\bpages{3519--3540}.
\bid{issn={1532-4435}, mr={2756192}}
\end{barticle}
%
\iffalse\OrigBibText
Ye, F. and Zhang, C. H. (2010).
Rate Minimaxity of the Lasso and
Dantzig Selector for the $\ell_q$ Loss in $\ell_r$ Balls.
\textit{J. Mach. Learn. Res.} \textbf{11}, 3519--3540.
\endOrigBibText\fi
\bptok{imsref}%
% NOT OUTPUTTED:
%   fjournal = Journal of Machine Learning Research (JMLR)
\endbibitem

%b29 ###
%b29 #&#
\bibitem{Yuanetal2009}
\begin{barticle}[mr]
\bauthor{\bsnm{Yuan},~\bfnm{Ming}\binits{M.}},
\bauthor{\bsnm{Joseph},~\bfnm{V.~Roshan}\binits{V.R.}} \AND
\bauthor{\bsnm{Zou},~\bfnm{Hui}\binits{H.}}
(\byear{2009}).
\btitle{Structured variable selection and estimation}.
\bjournal{Ann. Appl. Stat.}
\bvolume{3}
\bpages{1738--1757}.
\bid{doi={10.1214/09-AOAS254}, issn={1932-6157}, mr={2752156}}
\end{barticle}
%
\iffalse\OrigBibText
Yuan, M., Joseph, R. and Zou, H. (2009). Structured variable selection
and estimation. \textit{Ann. Appl. Stat.} \textbf{3}, 1738--1757.
\endOrigBibText\fi
\bptok{imsref}%
% NOT OUTPUTTED:
%   number = 4
%   doi = http://dx.doi.org/10.1214/09-AOAS254
%   fjournal = The Annals of Applied Statistics
\endbibitem

%b30 ###
%b30 #&#
\bibitem{YuanLin2006}
\begin{barticle}[mr]
\bauthor{\bsnm{Yuan},~\bfnm{Ming}\binits{M.}} \AND
\bauthor{\bsnm{Lin},~\bfnm{Yi}\binits{Y.}}
(\byear{2006}).
\btitle{Model selection and estimation in regression with grouped variables}.
\bjournal{J.~R.~Stat. Soc. Ser. B Stat. Methodol.}
\bvolume{68}
\bpages{49--67}.
\bid{doi={10.1111/j.1467-9868.2005.00532.x}, issn={1369-7412}, mr={2212574}}
\end{barticle}
%
\iffalse\OrigBibText
Yuan, M. and Lin, Y. (2006).
Model selection and estimation in regression with grouped variables.
\textit{ J. Roy. Stat. Soc. B. Met.} \textbf{68}, 49--67.
\endOrigBibText\fi
\bptok{imsref}%
% NOT OUTPUTTED:
%   number = 1
%   doi = http://dx.doi.org/10.1111/j.1467-9868.2005.00532.x
%   fjournal = Journal of the Royal Statistical Society. Series B. Statistical Methodology
\endbibitem

%b31 ###
%b31 #&#
\bibitem{ZhangHuang2008}
\begin{barticle}[mr]
\bauthor{\bsnm{Zhang},~\bfnm{Cun-Hui}\binits{C.-H.}} \AND
\bauthor{\bsnm{Huang},~\bfnm{Jian}\binits{J.}}
(\byear{2008}).
\btitle{The sparsity and bias of the LASSO selection in high-dimensional linear regression}.
\bjournal{Ann. Statist.}
\bvolume{36}
\bpages{1567--1594}.
\bid{doi={10.1214/07-AOS520}, issn={0090-5364}, mr={2435448}}
\end{barticle}
%
\iffalse\OrigBibText
Zhang, C. H. and Huang, J. (2008). The sparsity and bias of the Lasso
selection in high-dimensional linear regression. \textit{Ann.
Statist.} \textbf{36}, 1567--1594.
\endOrigBibText\fi
\bptok{imsref}%
% NOT OUTPUTTED:
%   number = 4
%   doi = http://dx.doi.org/10.1214/07-AOS520
%   coden = ASTSC7
%   fjournal = The Annals of Statistics
\endbibitem

%b32 ###
%b32 #&#
\bibitem{Zhaoetal2009}
\begin{barticle}[mr]
\bauthor{\bsnm{Zhao},~\bfnm{Peng}\binits{P.}},
\bauthor{\bsnm{Rocha},~\bfnm{Guilherme}\binits{G.}} \AND
\bauthor{\bsnm{Yu},~\bfnm{Bin}\binits{B.}}
(\byear{2009}).
\btitle{The composite absolute penalties family for grouped and hierarchical variable selection}.
\bjournal{Ann. Statist.}
\bvolume{37}
\bpages{3468--3497}.
\bid{doi={10.1214/07-AOS584}, issn={0090-5364}, mr={2549566}}
\end{barticle}
%
\iffalse\OrigBibText
Zhao, P., Rocha, G., and Yu, B. (2009).
The composite absolute penalties family for grouped and hierarchical
variable selection. \textit{Ann. Statist.} \textbf{37}, 3468--3497.
\endOrigBibText\fi
\bptok{imsref}%
% NOT OUTPUTTED:
%   number = 6A
%   doi = http://dx.doi.org/10.1214/07-AOS584
%   fjournal = The Annals of Statistics
\endbibitem

%b33 ###
%b33 #&#
\bibitem{ZhaoYu2006}
\begin{barticle}[mr]
\bauthor{\bsnm{Zhao},~\bfnm{Peng}\binits{P.}} \AND
\bauthor{\bsnm{Yu},~\bfnm{Bin}\binits{B.}}
(\byear{2006}).
\btitle{On model selection consistency of {L}asso}.
\bjournal{J. Mach. Learn. Res.}
\bvolume{7}
\bpages{2541--2563}.
\bid{issn={1532-4435}, mr={2274449}}
\bptnote{check pages}%
\end{barticle}
%
\iffalse\OrigBibText
Zhao, P. and Yu, B. (2006). On model selection consistency of Lasso.
\textit{J. Mach. Learn. Res.} \textbf{7}, 2541--2567.
\endOrigBibText\fi
\bptok{imsref}%
% NOT OUTPUTTED:
%   fjournal = Journal of Machine Learning Research (JMLR)
\endbibitem
\end{thebibliography}
\end{document}